\documentclass[twoside,11pt,reqno]{amsart}
\usepackage{a4wide}
\usepackage{hyperref}
\usepackage[dvips]{graphicx}
\usepackage{amsmath}
\usepackage[numbers]{natbib}
\usepackage{color}
\usepackage{amsmath,amssymb,bbm,enumerate}
\usepackage{pstricks}
\newtheorem{theorem}{Theorem}[section]

\newtheorem{corollary}{Corollary}[section]
\newtheorem{definition}{Definition}[section]
\newtheorem{lemma}{Lemma}[section]
\newtheorem{proposition}[theorem]{Proposition}
\newtheorem{remark}{Remark}[section]

\def\rmd{\mathrm{d}}
\def\rme{\mathrm{e}}
\def\rmi{\mathrm{i}}
\def\1{\mathbbm{1}}

\def\K{\mathbf{L}}

\def\bh{\mathbf{h}}
\def\bW{\mathbf{W}}
\def\bS{\mathbf{S}}

\def\bg{\mathbf{g}}
\def\bof{\mathbf{f}}
\def\bI{\mathbf{I}}
\def\bkappa{\boldsymbol{\kappa}}

\numberwithin{equation}{section}

\begin{document}

\title[Wavelet estimation of the long memory parameter]{Wavelet estimation of the long memory parameter for Hermite polynomial
of Gaussian processes}

\date{\today}

%    author one information
\author{M. Clausel}
\address{Laboratoire Jean Kuntzmann\\
Universit\'e de Grenoble, CNRS\\
F38041 Grenoble Cedex 9}
%\curraddr{}
\email{marianne.clausel@imag.fr}
%\thanks{}

%    author two information
\author{F. Roueff}
\address{Institut Mines--Telecom, Telecom ParisTech, CNRS LTCI, 46 rue Barrault\\
  75634 Paris Cedex 13, France}
%\curraddr{}
\email{roueff@telecom-paristech.fr}
\thanks{F. Roueff's research was
partially supported by the ANR project \emph{MATAIM} NT09 441552}

\author{M.~S. Taqqu}
\address{Departement of Mathematics and Statistics, Boston University,
  Boston, MA 02215, USA}
%\curraddr{}
\email{murad@math.bu.edu} \thanks{Murad~S.Taqqu was supported in part by the
  NSF grants DMS--0608669 and DMS--1007616 at Boston University.}

\author{C. Tudor}
\address{Laboratoire Paul Painlev\'e, UMR 8524 du CNRS, Universit\'e
  Lille 1, 59655 Villeneuve d'Ascq, France. Associate member: SAMM,
  Universit{\'e} de Panth\'eon-Sorbonne Paris 1.}
%\curraddr{}
\email{Ciprian.Tudor@math.univ-lille1.fr}
\thanks{C. Tudor's research was
  partially supported by the ANR grant \emph{Masterie} BLAN 012103.}

\subjclass[2010]{Primary  42C40 ; 60G18 ; 62M15 ; Secondary: 60G20, 60G22}

\keywords{Hermite processes ; Wavelet coefficients ; Wiener chaos ;
self-similar processes ; Long--range dependence.}

\begin{abstract}
  We consider stationary processes with long memory which are non--Gaussian and
  represented as Hermite polynomials of a Gaussian process. We focus on the
  corresponding wavelet coefficients and study the asymptotic behavior of the
  sum of their squares since this sum is often used for estimating the
  long--memory parameter.  We show that the limit is not Gaussian but can be
  expressed using the non--Gaussian Rosenblatt process defined as a
  Wiener-It\^{o} integral of order $2$. This happens even if the original
  process is defined through a Hermite polynomial of order higher than $2$.
\end{abstract}

\maketitle

\section{Introduction}\label{s:intro}

Wavelet analysis is a popular method for estimating the memory parameter of
stochastic processes with long--range dependence. The idea of using wavelets to
estimate the memory parameter $d$ goes back to~\cite{wornell:oppenheim:1992}
and~\cite{flandrin:1989a,flandrin:1989b,flandrin:1991a,flandrin:1999}.  See
also~\cite{abry:veitch:1998,abry:veitch:flandrin:1998}, \cite{bardet:2002},
\cite{bardet:bibi:jouini:2008},
\cite{bardet:lang:moulines:soulier:2000}. Wavelet methods are an alternative to
the Fourier methods developed by Fox and Taqqu~(\cite{fox:taqqu:1986}) and
Robinson~(\cite{robinson:1995:GPH,robinson:1995:GSE}).  The case of the
Gaussian processes, especially the fractional Brownian motion has been widely
studied. In this paper we will make an analysis of the wavelet coefficients of
stationary processes with long memory which are not Gaussian. The need for
non-Gaussian self-similar processes in practice (for example in hydrology) is
already mentioned in~\cite{Taq78} based on the study of stochastic modeling for
river-flow time series in~\cite{LK77}. More recently such an approach was used
for modeling Internet traffic, see~\cite[Chapter~3 and~4]{scherrer06:these}.

The wavelet analysis of non-Gaussian stochastic processes has been much less
treated in the literature. See~\cite{abry-helgason-pipiras-2011} for some
empirical studies.  \citeauthor{bardet:tudor:2010}, in
\cite{bardet:tudor:2010}, considered the case of the Rosenblatt process which
is a non-Gaussian self-similar process with stationary increments living in the
second Wiener chaos, that is, it can be expressed as a double iterated integral
with respect to the Wiener process.  It can be also defined as a Hermite
process of order 2, while the fractional Brownian motion is a Hermite process
of order 1. We refer to Section \ref{mainresult} for the definition of the
Rosenblatt process (see also \cite{Taq79}, \cite{abry:pipiras:2006},
\cite{GS85}), and to \cite{breton:nourdin:2008}, \cite{CTV10},
\cite{embrechts:maejima:2002}, \cite{Taq79} for the definition and various
properties of the Hermite process.

% We consider processes with long memory that are not Gaussian and study the
% asymptotic behavior of their scalogram. The scalogram is a sum of squares of
% wavelet coefficients which is used to estimate the intensity of long memory.
% A study of Gaussian and linear processes with long memory were developed in
% \cite{moulines:roueff:taqqu:2007:jtsa,roueff-taqqu-2009,roueff-taqqu-2009b}.
% In this paper, we consider non-Gaussian processes defined as Hermite polynomials
% of arbitrary order of a Gaussian process.  It is natural to focus on Hermite
% polynomials since these forms an orthogonal basis of the $L^2$ space generated
% by a Gaussian process.

% \cite{bardet:tudor:2010} considered the case of the Rosenblatt process which is
% a non-Gaussian self-similar process with stationary increments.
% Recall that $\{X_t,\,t\in\mathbb{R}\}$ is self-similar with index $H$ ICI
% leaving in the second Wiener chaos, that is, it can be expressed
% as a double iterated integral with respect to the Wiener process.
% It can be also defined as a Hermite process of order 2, while the
% fractional Brownian motion is a Hermite process of order 1. We
% refer to Section \ref{mainresult} for the definition  of the
% Rosenblatt process (see also \cite{Taq79}, \cite{Tud08}), and to
% \cite{CTV10}, \cite{pipiras:taqqu2010}, \cite{Taq79}  for the
% definition and various properties of the Hermite process.

In the present work, we consider processes expressed as a Hermite polynomial of
order greater than $1$ of a Gaussian time series.  This will allow us to gain
insight into more complicated situations. A more general case, involving
processes that can be expressed as (finite or infinite) sum of Hermite
polynomials of a Gaussian time series is studied in our recent work
\cite{clausel-roueff-taqqu-tudor-2012}.  In
this work, we use :
\begin{enumerate}[a)]
\item a wide class of wavelets as in~(\ref{e:W1}), instead of ``variations'';
\item an input process with long-range dependence, as in~(\ref{e:dq>0}), instead
  of self-similar processes;
\item a semiparametric setup, as in~(\ref{e:sdf}), instead of a parametric one.
\end{enumerate}
We derive the limit theorems that are needed for wavelet--based estimation
procedures of the memory parameter.  We will investigate the estimation problem
in another paper.

Denote by $X=\{X_{t}\}_{t\in\mathbb{Z}}$ a centered stationary Gaussian process
with unit variance and spectral density $f(\lambda), \lambda \in (-\pi ,
\pi)$. Such a stochastic process is said to have {\it short memory} or {\it
  short--range dependence} if $f(\lambda)$ is positive and bounded around
$\lambda=0$ and {\it long memory} or {\it long--range dependence} if
$f(\lambda)\to\infty$ as $\lambda\to 0$. We will suppose that
$\{X_{t}\}_{t\in\mathbb{Z}}$ has long--memory with memory parameter $0<d<1/2$,
that is,
\begin{equation}\label{e:sdf0}
f(\lambda)\sim |\lambda|^{-2d}f^*(\lambda)\mbox{ as }\lambda \to 0
\end{equation}
where $f^*(\lambda)$ is a bounded spectral density which is continuous and
positive at the origin.  It is convenient to set
\begin{equation}\label{e:sdf}
f(\lambda)=|1-\rme^{-\rmi\lambda}|^{-2d}f^*(\lambda),\quad\lambda\in (-\pi,\pi]\;.
\end{equation}
Since the spectral density of a stationary process is integrable, we require
$d<\frac{1}{2}$.

We shall also consider a process $\{Y_{t}\}_{t\in\mathbb{Z}}$, not necessarily
stationary but its difference $\Delta^K Y$ of order $K\geq0$ is
stationary. Moreover, instead of supposing that $\Delta^K Y$ is Gaussian, we
will assume that
\begin{equation}\label{e:K}
\left(\Delta^K Y\right)_{t}=H_{q_0}(X_{t}),\quad t\in\mathbb{Z}\;,
\end{equation}
where $(\Delta Y)_{t}=Y_{t}-Y_{t-1}$, where $X$ is
Gaussian with spectral density $f$ satisfying~(\ref{e:sdf}) and where $H_{q_0}$ is the $q_0$--th Hermite polynomial.

We will focus on the wavelet coefficients of
$Y=\{Y_{t}\}_{t\in\mathbb{Z}}$. Since $\{Y_{t}\}_{t\in\mathbb{Z}}$
is random so will be its wavelet coefficients which we denote by
$\{W_{j,k},\,j\geq 0,\,k\in\mathbb{Z}\}$, where $j$ indicates the
scale and $k$ the location. These wavelet coefficients are defined
by
\begin{equation}
  \label{eq:wav_coeff_def1}
 W_{j,k}=\sum_{t\in\mathbb{Z}}h_j(\gamma_j k-t)Y_{t}\;,
\end{equation}
where $\gamma_j\uparrow \infty$ as $j\uparrow \infty$ is a
sequence of non--negative scale factors applied at scale $j$, for
example $\gamma_j=2^j$ and $h_{j}$ is a filter whose properties
will be listed below. We follow the engineering convention where
large values of $j$ correspond to large scales. Our goal is to
find the distribution of the empirical quadratic mean of these
wavelet coefficients at large scales $j\to\infty$, that is, the
asymptotic behavior of the scalogram
\begin{equation}\label{e:defsnj}
S_{n,j}=\frac{1}{n}\sum_{k=0}^{n-1}W_{j,k}^2\;,
\end{equation}
adequately normalized as the number of wavelet coefficients $n$ and
$j=j(n)\to\infty$. This is a necessary and important step in developing methods
for estimating the underlying long memory parameter $d$, see the references
mentioned at the beginning of this section. Indeed, using the wavelet
scalogram, there is standard way to construct an estimator of the memory
parameter. The asymptotic behavior of the scalogram gives the convergence rate
of this estimator. We provide more details in Section
\ref{sec:estim-long-memory}.

When $q_0=1$, the behavior of $S_{n,j}$ has been studied
in~\cite{roueff-taqqu-2009b}. In this case, under certain conditions, the limit
as $j,n\to \infty$ of the suitably renormalized sequence $S_{n,j}$ is
Gaussian. If $q_0\geq 2$ only few facts are known on the behavior of the
scalogram $S_{n,j}$.  In~\cite{bardet:tudor:2010}, the authors have made a
wavelet analysis of the Rosenblatt process (see Definition \ref{d:HP} with
$q=2$). This situation roughly corresponds to the case $q_0=2$ (the second
Hermite polynomial). It has been shown that its associated scalogram has a
non-Gaussian behavior, that is, after normalization it converges to a
Rosenblatt random variable.  Basically, what happens is the following: the
random variable $H_{2}(X_t)$ is, for every $t \in \mathbb{Z}$ an element of the
second Wiener chaos and its square can be decomposed, using the properties of
multiple stochastic integrals, as a sum of a multiple integral in the fourth
Wiener chaos and a multiple integral in the second Wiener chaos.  It turns out
that the leading term is the one in the second Wiener chaos which converges to
a Rosenblatt random variable (a Rosenblatt process at time 1).  Wavelet
analysis for $G=H_{q}$ with $q>2$ has not been done until now. Some intuition
can be gained from the study of quadratic variations of the increments of the
Hermite process, in~\cite{CTV10}. In this case the starting process is
self--similar, that is, invariant under scaling. Again the limit turns out to
be the Rosenblatt random variable. Briefly since the Hermite process is an
element of the $q$th Wiener chaos, its square (minus the expectation of its
square) can be expressed as a sum of multiple integrals of orders 2,4,.. until
$2q$. It turns out that the main term is the one in the second Wiener chaos
which converges to a Rosenblatt random variable. This may suggest that in our
situation one would have perhaps a ``reduction theorem'' as
in~\cite{dobrushin:major:1979}, stating that it is the lower order term which
dominates. This is not the case however. We will show in a subsequent paper
that higher--order Hermite processes can appear in the limit even when the
initial data are a mixture of a Gaussian and non--Gaussian components. See also
\cite{nourdin:peccati:2010}, \cite{nourdin:peccati:2009} for other examples of
limit theorems based on the chaos expansion.

The paper is structured as follows. In Section~\ref{s:wavcoeff} we introduce
the wavelet filters and state the assumptions imposed on them. In
Section~\ref{mainresult} we state our main result and we introduce the
Rosenblatt process which appears as limit for $q_0\geq 2$.  This result is
stated for a multivariate scalogram considered at a single scale.  In
Section~\ref{sec:from-mult-scal}, we explain how this applies to the asymptotic
behavior of the univariate scalogram at multiple scales (in short, the
\emph{multiscale asymptotics}). Results on the estimation of the long memory
parameter are derived in Section~\ref{sec:estim-long-memory}.  In
Section~\ref{s:expansion} we give the chaos expansion of the
scalogram. Section~\ref{s:L2bound} and~\ref{s:leadingterms} describe the
asymptotic behavior of the main terms appearing in the decomposition of the
scalogram. The proof of the main results is in
Section~\ref{s:proofs-main}.  Finally, Sections~\ref{sec:techlemma} contains
technical lemmas used throughout our paper and Appendix~\ref{s:appendix}
recalls the basic facts needed in this paper about Wiener chaos.

\section{The wavelet coefficients}\label{s:wavcoeff}

The Gaussian sequence $X=\{X_t\}_{t\in\mathbb{Z}}$ with spectral
density~(\ref{e:sdf}) is long--range dependent because $d>0$ and
hence its spectrum explodes at $\lambda=0$. Whether
$\{H_{q_0}(X_t)\}_{t\in\mathbb{Z}}$ is also long-range dependent
depends on the respective values of $q_0$ and $d$.  We show
in~\cite{clausel-roueff-taqqu-tudor-2011a}, that the spectral density of
$\{H_{q_0}(X_t)\}_{t\in\mathbb{Z}}$ behaves proportionally to
$|\lambda|^{-\delta_+(q_0)}$ as $\lambda\to 0$, where
\begin{equation}\label{e:ldparamq}
\delta_+(q)=\max(\delta(q),0)\quad\text{and}
\quad\delta(q)=q d-(q-1)/2,\qquad q=1,2,3,\dots\;,
\end{equation}
and hence $\delta_+(q_0)$ is the memory parameter of
$\{H_{q_0}(X_t)\}_{t\in\mathbb{Z}}$ . Therefore, since $0<d<1/2$,
in order for $\{H_{q_0}(X_t)\}_{t\in\mathbb{Z}}$, $q_0\geq 1$, to
be long--range dependent, one needs
\begin{equation}\label{e:dq>0}
\delta(q_0)>0\Leftrightarrow (1-1/q_0)/2<d<1/2\;,
%\quad\mbox{ that is if and only if }q<1/(1-2d)\;.
\end{equation}
that is, $d$ must be sufficiently close to $1/2$. Specifically,
for long--range dependence,
\[
q_0=1\Rightarrow d>0,\quad q_0=2\Rightarrow d>1/4,\quad q_0=3\Rightarrow d>1/3,\quad q_0=4\Rightarrow d>3/8\;\dots
\]
From another perspective, for all $q_0\geq1$
\begin{equation}\label{e:dqq>0}
\delta(q_0)>0\Leftrightarrow q_0<1/(1-2d)\;,
\end{equation}
and thus $\{H_{q_0}(X_t)\}_{t\in\mathbb{Z}}$ is short--range dependent
if $q_0\geq 1/(1-2d)$.
In the following, we always assume
that $\{H_{q_0}(X_t)\}_{t\in\mathbb{Z}}$ has
long memory, that is,
\begin{equation}\label{e:longmemorycondition}
1\leq q_0 < 1/(1-2d)\text{ or, equivalently, }0<\delta(q_0)<1/2 \;.
\end{equation}

As indicated in the introduction, we consider the process
$\{Y_t\}_{t\in\mathbb{Z}}$, where $\Delta^K Y_t=H_{q_0}(X_t)$ for
any $t\in\mathbb{Z}$ and for some $K\geq 0$ (see~(\ref{e:K})). We
are interested in the wavelets coefficients of the process
$\{H_{q_0}(X_t)\}_{t\in\mathbb{Z}}$. To obtain them, one applies a
linear filter $h_j(\tau),\tau\in\mathbb{Z}$, at each scale $j\geq
0$. We shall characterize below the filters $h_j(\tau)$ by their
discrete Fourier transform~:
\begin{equation}\label{e:dF}
\widehat{h}_j(\lambda)=\sum_{\tau\in\mathbb{Z}}h_j(\tau)
\rme^{-\rmi \lambda\tau },\,\lambda\in [-\pi,\pi]\;,\quad h_j(\tau)=\frac{1}{2\pi}\int_{-\pi}^{\pi}\widehat{h}_j(\lambda)\rme^{\rmi\lambda\tau}\rmd\lambda,\tau\in\mathbb{Z}\;.
\end{equation}
The resulting wavelet coefficients $W_{j,k}$, where $j$ is the scale and $k$ the location are defined as
\begin{equation}\label{e:W1}
W_{j,k}=\sum_{t\in\mathbb{Z}}h_j(\gamma_j k-t)Y_{t}=\sum_{t\in\mathbb{Z}}h_j(\gamma_j k-t)\Delta^{-K}H_{q_0}(X_{t}),\,j\geq 0, k\in\mathbb{Z},
\end{equation}
where $\gamma_j\uparrow \infty$ as $j\uparrow \infty$ is a sequence of
non--negative scale factors applied at scale $j$, for example $\gamma_j=2^j$.
We do not assume that the wavelet coefficients are orthogonal nor that they are
generated by a multiresolution analysis, but only that the filters $h_j$
concentrate around the zero frequency as $j\to\infty$ with some uniformity, see
Assumptions (W-\ref{ass:w-b})--(W-\ref{ass:w-c}) below.

To study the joint convergence at several scales jointly going to infinity,
wavelet coefficients can be considered as a process $W_{j+m_0,k}$ indexed by
$m_0,k$ and where we let $j\to\infty$ as in
\cite{clausel-roueff-taqqu-tudor-2011a}. Here we are interested in the
scalogram defined as the empirical square mean~(\ref{e:defsnj}) with $n$ equal
to the number of wavelets coefficients at scale $j$ available from $N$
observations of the original process $Y_1,\dots,Y_N$. Considering the joint
asymptotic behavior at various scales means that we have to deal with different
down-sampling rates $\gamma_j$ and different numbers $n_j$ of available wavelet
coefficients, both indexed by the scale $j$. It is shown in
\cite{roueff-taqqu-2009b} that the joint behavior of the scalogram at multiple
scales can be deduced from the joint behavior of the
statistic~(\ref{e:defsnj}), viewed as a vector whose components have the same
$j$ and $n$ but different filters $h_{\ell,j}$, $\ell=1,\dots,m$. We shall
adopt the multivariate scalogram setup in our asymptotic analysis. We shall
apply it
in Section~\ref{sec:from-mult-scal} to deduce the
multiscale asymptotic behavior of the univariate scalogram. This will also
allow us to contrast the cases $q_0>1$ treated in this contribution with the
case $q_0=1$ which follows from the result obtained in
\cite{roueff-taqqu-2009}.
% \textcolor{red}{SUPPRIMER To illustrate how this is done, suppose, for
%   example, that we want to study the joint behavior of $W_{j+m_0,k}$ for
%   $m_0\in\{0,-1\}$. Recall that $j-1$ is a finer scale than $j$. Following the
%   framework of \cite{roueff-taqqu-2009b}, one could consider the multivariate
%   coefficients $\bW_{j,k}=(W_{j,k},\,W_{j-1,2k},\,W_{j-1,2k+1})$ since there
%   are twice as many wavelet coefficients at scale $j-1$ as at scale $j$.  These
%   coefficients can be considered as the output of a multidimensional filter
%   $\bh_j$ defined as
%   $\bh_j(\tau)=(h_j(\tau),\,h_{j-1}(2\tau),\,h_{j-1}(2\tau+2^{j-1}))$.}  \textcolor{red}{SUPPRIMER, and stated in the special
%   case of the scalogram in \cite[Theorem~1]{roueff-taqqu-2009b}.}
Our
assumption on the filters $h_{\ell,j}$, $\ell=1,\dots,m$ are the same as in \cite[Theorem~1]{roueff-taqqu-2009b}, except
that we allow $\gamma_j\neq2^j$ for the sake of generality, and we assume
locally uniform convergence in the asymptotic behavior in~(\ref{EqLimHj}).
These assumptions are satisfied in the standard wavelet analysis described in
\cite{moulines:roueff:taqqu:2007:jtsa} and   briefly referred to in
Section~\ref{sec:from-mult-scal}.

From now on, the wavelet coefficient $W_{j,k}$ defined in~(\ref{e:W1}) will be
supposed to be $\mathbb{R}^m$-valued with $h_j$ representing a $m$-dimensional
vector with entries $h_{\ell,j}$, $\ell=1,\dots,m$. We will use bold faced
symbols $\bW_{j,k}$ and $\bh_j$ to emphasize the multivariate setting, thus
\begin{equation}\label{e:bW1}
\bW_{j,k}=\sum_{t\in\mathbb{Z}}\bh_j(\gamma_j
k-t)Y_{t}=\sum_{t\in\mathbb{Z}}\bh_j(\gamma_j
k-t)\Delta^{-K}H_{q_0}(X_{t}),\,j\geq 0, k\in\mathbb{Z}.
\end{equation}
  We shall make the following assumptions on the filters $\bh_j$:
\begin{enumerate}[(W-a)]
\item\label{ass:w-a} \underline{Finite support}: For each $\ell$ and $j$,
  $\{h_{\ell,j}(\tau)\}_{\tau\in\mathbb{Z}}$ has finite support.
\item\label{ass:w-b} \underline{Uniform smoothness}: There exists $M\geq 0$,
  $\alpha>1/2$ and $C>0$ such that for all $j\geq0$ and
  $\lambda\in [-\pi,\pi]$,
    \begin{equation}\label{e:majoHj}
    |\widehat{\bh}_{j}(\lambda)|\leq \frac{C\gamma_j^{1/2}|\gamma_j\lambda|^M}{(1+\gamma_j|\lambda|)^{\alpha +M}}\;,
    \end{equation}
    where $|x|$ denotes the Euclidean norm of vector $x$. By $2\pi$-periodicity
    of $\widehat{\bh}_{j}$ this inequality can be extended to
    $\lambda\in\mathbb{R}$ as
\begin{equation}\label{EqMajoHjR}
|\widehat{\bh}_{j}(\lambda)|\leq C
\frac{\gamma_j^{1/2}|\gamma_j\{\lambda\}|^M}
{(1+\gamma_j|\{\lambda\}|)^{\alpha+M}}\;,
\end{equation}
where $\{\lambda\}$ denotes the element of $(-\pi,\pi]$
such that $\lambda-\{\lambda\}\in2\pi\mathbb{Z}$.
\item\label{ass:w-c} \underline{Asymptotic behavior}: There exist a sequence of
  phase
  functions $\Phi_j:\mathbb{R}\to(-\pi,\pi]$ and some function
  $\widehat{\bh}_{\infty}:\mathbb{R}\to\mathbb{C}^p$  such that
\begin{equation}\label{EqLimHj}
\lim_{j \to +\infty}\gamma_j^{-1/2}\widehat{\bh}_{j}(\gamma_j^{-1}\lambda)
\rme^{\rmi\Phi_j(\lambda)}
= \widehat{\bh}_{\infty}(\lambda)\;,
\end{equation}
 locally uniformly on $\lambda\in\mathbb{R}$.
\end{enumerate}
In (W--\ref{ass:w-c}), \emph{locally uniformly} means that
for all $r>0$,
$$
\sup_{|\lambda|\leq r}\left|\gamma_j^{-1/2}\widehat{\bh}_{j}(\gamma_j^{-1}\lambda)
\rme^{\rmi\Phi_j(\lambda)}- \widehat{\bh}_{\infty}(\lambda)\right|\to 0\;.
$$
This is satisfied if the set of filters correspond to a discrete wavelet
transform (see Proposition~3 in \cite{moulines:roueff:taqqu:2007:jtsa}).
Assumptions (\ref{e:majoHj}) and~(\ref{EqLimHj}) imply that for any
$\lambda\in\mathbb{R}$,
\begin{equation}\label{EqMajoHinf}
|\widehat{\bh}_{\infty}(\lambda)|\leq C\frac{|\lambda|^M}{(1+|\lambda|)^{\alpha+M}}\;.
\end{equation}
Hence vector $\widehat{\bh}_{\infty}$ has entries in $L^{2}(\mathbb{R})$. We
let  $\bh_{\infty}$ be the vector of $L^2(\mathbb{R})$ inverse Fourier
transforms of
$\widehat{h}_{\ell,\infty}$, $\ell=1,\dots,m$, that is
\begin{equation}\label{eq:TFdef}
\widehat{\bh}_{\infty}(\xi)=
\int_{\mathbb{R}} \bh_{\infty}(t)\rme^{-\rmi t\xi}\;\rmd t,
\quad \xi\in\mathbb{R}\;.
\end{equation}

Observe that while $\widehat{\bh}_{j}$ is $2\pi$--periodic, the function
$\widehat{\bh}_{\infty}$ has non--periodic entries on $\mathbb{R}$. For the
connection between these assumptions on $\bh_{j}$ and corresponding assumptions
on the scaling function $\varphi$ and the mother wavelet $\psi$ in the
classical wavelet setting see \cite{moulines:roueff:taqqu:2007:jtsa} and
\cite{roueff-taqqu-2009b}. In particular, in the univariate setting $m=1$, one
has $\widehat{h}_{\infty}=\widehat{\varphi}(0)\overline{\widehat{\psi}}$.

For $M\geq K$, a more convenient way to express $\bW_{j,k}$ is to incorporate
the linear filter $\Delta^{-K}$ in~(\ref{e:bW1}) into the filter $\bh_j$ and
denote the resulting filter $\bh_j^{(K)}$. Then
\begin{equation}\label{e:W2}
\bW_{j,k}=\sum_{t\in\mathbb{Z}}\bh_j^{(K)}(\gamma_j k-t)H_{q_0}(X_{t})\;,
\end{equation}
where
\begin{equation}\label{e:HjkVSHj}
\widehat{\bh}_j^{(K)}(\lambda)=(1-\rme^{-{\rmi}\lambda})^{-K} \widehat{\bh}_j(\lambda)
\end{equation}
is the component wise discrete Fourier transform of $\bh_j^{(K)}$.
Since $\{H_{q_0}(X_t),t\in\mathbb{Z}\}$ is stationary, so is
$\{W_{j,k},k\in\mathbb{Z}\}$ for each scale $j$.
Using~(\ref{EqMajoHjR}), we further get,
\begin{equation}\label{e:hK}
\left| \widehat{\bh}_j^{(K)}(\lambda) \right| \leq C
  \gamma_j^{1/2+K}\;\frac{|\gamma_j\{\lambda\}|^{M-K}}{(1+\gamma_j|\{\lambda\}|)^{\alpha+M}},\quad\lambda\in\mathbb{R},\,j\geq 1\;.
\end{equation}
In particular, if $M\geq K$, using that
$(|\gamma_j\{\lambda\}|/(1+\gamma_j|\{\lambda\}|))^{M}
\leq(|\gamma_j\{\lambda\}|/(1+\gamma_j|\{\lambda\}|))^{K}$,
we get
\begin{equation}
  \label{eq:M_replaced_by_K}
  \left| \widehat{\bh}_j^{(K)}(\lambda) \right| \leq C
  \gamma_j^{1/2+K}\;(1+\gamma_j|\{\lambda\}|)^{-\alpha-K},
  \quad\lambda\in\mathbb{R},\,j\geq 1\;.
\end{equation}

By Assumption~(\ref{e:majoHj}), $\bh_j$ has vanishing moments up to order
$M-1$, that is, for any integer  $0\leq k\leq M-1$,
\begin{equation}\label{e:mom}
\sum_{t\in \mathbb{Z}}\bh_j(t)t^k =0\;.
\end{equation}
Observe that $\Delta^K Y$ is centered by definition. However, by~(\ref{e:mom}),
the definition of $\bW_{j,k}$ only depends on $\Delta^M Y$. In
particular, provided that $M\geq K+1$, its value is not modified if a constant
is added to $\Delta^K Y$, whenever $M\geq K+1$.

\section{Main result}\label{mainresult}

Recall that
$$
(\Delta^KY)_t=H_{q_0}(X_t),\quad t\in\mathbb{Z}\;.
$$
The condition~(\ref{e:longmemorycondition}) ensures such that
$\{H_{q_0}(X_t)\}_{t\in\mathbb{Z}}$ is long-range dependent
(see~\cite{clausel-roueff-taqqu-tudor-2011a}, Lemma 4.1). Our main result deals
with the asymptotic behavior of the scalogram $S_{n,j}$, defined in the
univariate case $m=1$ by~(\ref{e:defsnj}) as $j,n\to\infty$, that is, as
$n\to\infty$ (large sample behavior) with $j=j(n)$ being an arbitrary diverging
sequence (large scale behavior). More precisely, we will study the asymptotic
behavior of the sequence
\begin{equation}\label{e:snjm}
\overline{\bS}_{n,j}%=\bS_{n,j}-\mathbb{E}[\left| \bW_{j,0}\right|^{2}]
=\frac{1}{n}\sum_{k=0}^{n-1}\left(\bW_{j,k}^2-\mathbb{E}[\bW_{j, k}^{2}]\right)
=\left[\frac{1}{n}\sum_{k=0}^{n-1}\left(W_{\ell,j,k}^2-\mathbb{E}[W_{\ell,j, k}^{2}]\right)\right]_{\ell=1,\dots,m}\;,
\end{equation}
adequately normalized as $j,n\to\infty$, where $W_{\ell,j,k}$,
$\ell=1,\dots,m$, denote the $m$ entries of vector $\bW_{j,k}$.
The limit will be
expressed in terms of the Rosenblatt process which is defined as
follows.
\begin{definition}\label{d:HP}
The Rosenblatt process of index $d$ with
\begin{equation}\label{e:qd}
1/4<d<1/2\;,
\end{equation}
is the continuous time process
\begin{equation}\label{e:harmros}
Z_{d}(t)= \int_{\mathbb{R}^{2}}^{\prime\prime} \frac{\rme^{\rmi
(u_1+u_2)\,t}- 1} {\rmi(u_1+u_2)}|u_1|^{-d}|u_2|^{-d}\;
\rmd\widehat{W}(u_1) \rmd\widehat{W}(u_2),\,t\in\mathbb{R}\;.
\end{equation}
\end{definition}
The multiple integral~(\ref{e:harmros}) with respect to the complex-valued Gaussian random measure $\widehat{W}$ is defined in
Appendix~\ref{s:appendix}. The symbol
$\int_{\mathbb{R}^{2}}^{\prime\prime}$ indicates that one does not
integrate on the diagonal $u_1=u_2$. The integral is well-defined
when~(\ref{e:qd}) holds because then it has finite $L^2$ norm. This
process is self--similar with self-similarity parameter
\[
H=2d \in (1/2,1),
\]
that is for all $a>0$, $\{Z_{d}(at)\}_{t\in\mathbb{R}}$ and
$\{a^H Z_{d}(t)\}_{t\in\mathbb{R}}$ have the same finite
dimensional distributions, see~\cite{Taq79}.\\

We now list the assumptions behind our main result:\\

\medskip
\noindent{\bf Assumptions A} $\{\bW_{j,k},\,j\geq1,k\in\mathbb{Z}\}$ are the
wavelet coefficients defined by~(\ref{e:bW1}) , where
\begin{enumerate}[(i)]
\item\label{item:A1} $X$ is a stationary Gaussian process with spectral density $f$
  satisfying~(\ref{e:sdf}) with $0<d<1/2$;
\item\label{item:A2} $H_{q_0}$ is the $q_0$ th Hermite polynomial where $q_0$
satisfies condition~(\ref{e:longmemorycondition});
\item\label{item:A3} the sequence of positive integers $(\gamma_j)_{j\geq1}$ is non-decreasing
  and diverging;
\item\label{item:A4} the wavelet filters $\bh_j=[h_{\ell,j}]_{\ell=1,\dots,m}$,
  ${j\geq1}$, satisfy (W-\ref{ass:w-a})--(W-\ref{ass:w-c}).
\end{enumerate}
The definition of Hermite polynomials is recalled in Appendix~\ref{s:appendix}.
The following theorem gives the limit of~(\ref{e:snjm}), suitably normalized,
as the number of wavelet coefficients and the scale $j=j(n)$ tend to infinity,
in the cases $q_0=1$ and $q_0\geq2$.
\begin{theorem}\label{thm:LD}
  Suppose that
  Assumptions \textbf{A} hold with $M\geq K+\delta(q_0)$, where $\delta(\cdot)$
  is defined in~(\ref{e:ldparamq}). Define the centered
  multivariate scalogram $\overline{\bS}_{n,j}$ by~(\ref{e:snjm}) and let $(n_j)$ be any
  diverging sequence of integers.
  \begin{enumerate}[(a)]
  \item\label{item:thm-case-gauss} Suppose $q_0=1$ and that
    $(\gamma_j)$ is a sequence of even integers. Then, as $j\to\infty$,
\begin{equation}
\label{eq:CenteredZ}
n_j^{1/2}\gamma_j^{-2(d+K)}
\overline{\bS}_{n_j,j}
\overset{\mathcal{L}}{\longrightarrow}
\mathcal{N}(0,\Gamma) \; ,
\end{equation}
where $\Gamma$ is the $m\times m$ matrix with entries
\begin{equation}
  \label{eq:GammaDef}
  \Gamma_{\ell,\ell'}= 4\pi (f^*(0))^2 \;  \;
\int_{-\pi}^\pi \left|
\sum_{p\in\mathbb{Z}}
|\lambda+2p\pi|^{-2(K+d)}
[\widehat{h}_{\ell,\infty}\overline{\widehat{h}_{\ell',\infty}}](\lambda+2p\pi)
\right|^2 \, \rmd\lambda\;,\quad 1\leq \ell, \ell'\leq m\; .
\end{equation}
\item\label{item:thm-case-rozen} Suppose $q_0\geq 2$.
Then as $j\to\infty$,
\begin{equation}
\label{eq:caseBasymp_dist}
n_j^{1-2d} \gamma_{j}^{-2(\delta(q_{0})+K)} \overline{\bS}_{n_j,j}
\overset{\mathcal{L}}{\longrightarrow} f^*(0)^{q_0}\,\K_{q_0-1} \,
  Z_d(1) \;,
\end{equation}
where $Z_d(1)$ is the Rosenblatt process in~(\ref{e:harmros}) evaluated at time
$t=1$, $f^*(0)$ is the \emph{short-range spectral density} at zero frequency
in~(\ref{e:sdf0}) and where $\K_{q_0-1}$ is the deterministic $m$-dimensional
vector $[L_{q_0-1}(\widehat{h}_{\ell,\infty})]_{\ell=1,\dots,m}$ with finite
entries defined by
\begin{equation}\label{e:defKp}
L_{p}(g)=\int_{\mathbb{R}^{p}}\frac{|g(u_1+\cdots+u_{p})|^{2}}
{|u_1+\cdots+u_{p}|^{2K}}\; \prod_{i=1}^{p} |u_i|^{-2d}\; \rmd
u_1\cdots\rmd u_{p} \;,
\end{equation}
for any $g:\mathbb{R}\to\mathbb{C}$ and $p\geq1$.
  \end{enumerate}
\end{theorem}
This theorem is proved in Section~\ref{s:proofs-main}.
\begin{remark}
  Since $\delta(1)=d$ we observe that the exponent of $\gamma_j$ in the rate of
  convergence of $\overline{\bS}_{n,j}$ can be written as $-2(\delta(q_{0})+K)$
  for both cases $q_0=1$ and $q_0\geq2$, see~(\ref{eq:CenteredZ})
  and~(\ref{eq:caseBasymp_dist}), respectively. This corresponds to the fact
  that $d_0=\delta(q_0)+K$ is the long memory parameter of $Y$, and, as a
  consequence, $\mathbb{E}|\bW_{j,0}|^2\sim C \gamma_j^{2(\delta(q_0)+K)}$ as
  $j\to\infty$, see for example Theorem~5.1 in
  \cite{clausel-roueff-taqqu-tudor-2011a}.  In contrast, the exponent of $n$ is
  always larger in the case $q_0\geq2$, since this implies $2d-1>-1/2$ under
  Condition~(\ref{e:longmemorycondition}).  The statistical behavior of the
  limits are also very different in the two cases. In~(\ref{eq:CenteredZ}) the
  limit is Gaussian while in~(\ref{eq:caseBasymp_dist}), the limit is
  Rosenblatt. Another difference is that the entries of the limit vector
  in~(\ref{eq:caseBasymp_dist}) have cross-correlations equal to 1 (they only
  differ through a multiplicative constant). In contrast, this typically does
  not happen in~(\ref{eq:CenteredZ}). %By **Lemma???**
% the Cauchy-Schwartz Inequality, this would happen
%   only if
%   $\widehat{h}_{i,\infty}(\lambda+2p\pi)/\widehat{h}_{i',\infty}(\lambda+2p\pi)$
%   does not depend on $p\in\mathbb{Z}$ for almost every $\lambda$.
\end{remark}

\begin{remark}
While  $H_{q_0}(X_t)$ involves a single multiple integral
of order $q_0$, $\bW_{j,k}^2$ and hence $\overline{\bS}_{n,j}$
in~(\ref{e:snjm}) involves a sum of multiple integrals of order 0, 2, 4, 6...
up to $2q_0$. But the limiting Rosenblatt process in Theorem~\ref{thm:LD}
involves only a double integral, albeit with a non--random factor
$\K_{q_0-1}$ expressed as a non--random multiple integral of order
$q_0-1$. In view of Theorem~5.1 of \cite{clausel-roueff-taqqu-tudor-2011a}, the
components of $\K_{q_0-1}$ are the asymptotic variances of the wavelet
coefficients applied to $\Delta^{-K}H_{q_0-1}(X_t)$.
\end{remark}

\section{From multivariate to multiscale asymptotics}
\label{sec:from-mult-scal}

Theorem~\ref{thm:LD} applies to multivariate filters $\bh_j$
which define the scalogram $\bS_{n,j}$.  We will use it to obtain in
Theorem~\ref{thm:multiscale} multiscale asymptotics for univariate filters and
corresponding scalograms. This passage between these two prospectives is
explained in the proof of Theorem~\ref{thm:multiscale}. We use dyadic scales
here, as in the standard wavelet analysis described in
\cite{moulines:roueff:taqqu:2007:jtsa}, where the wavelet coefficients are
defined as
\begin{equation}
  \label{eq:wav_coeff_def1_dyad}
 W_{j,k}=\sum_{t\in\mathbb{Z}}g_j(2^j k-t)Y_{t}\;,
\end{equation}
which corresponds to~(\ref{eq:wav_coeff_def1}) with $\gamma_j=2^j$ and with
$(g_j)$ denoting a sequence of filters that
satisfies~(W-\ref{ass:w-a})--(W-\ref{ass:w-c}) with $m=1$, and $M$ and $\alpha$
respectively defined as the number of vanishing moments of the wavelet and its
Fourier decay exponent.
In the case of a multiresolution analysis, $g_j$ can be
deduced from the associated mirror filters.

The number $n_j$ of wavelet coefficients available at
scale $j$, is related both to the number $N$ of observations $Y_1,\cdots,Y_N$
of the time series $Y$ and to the length $T$ of the support of the wavelet
$\psi$. More precisely, one has
\begin{equation}
  \label{eq:nj_def}
  n_j=[2^{-j}(N-T+1)-T+1]= 2^{-j}N+0(1)\;,
\end{equation}
where $[x]$ denotes the integer part of $x$ for any real $x$.
Details about the above facts can be found in
\cite{moulines:roueff:taqqu:2007:jtsa,roueff-taqqu-2009b}.

In this context, the scalogram is an empirical measure of the distribution of
``energy of the
signal'' along scales, based on the $N$ observations $Y_1,\cdots,Y_N$. It is
defined as
\begin{equation}
  \label{eq:scalo_def}
\widehat{\sigma}^2_{j}=\frac1{n_j}\sum_{k=0}^{n_j-1}W_{j,k}^2, \quad j\geq0 \;,
\end{equation}
and is identical to $S_{n_j,j}$ defined in~(\ref{e:defsnj}).
Note that the sequence $(\widehat{\sigma}^2_{j})_{j\geq0}$ is indexed by the
scale index $j$ but also depends on the number $N$ of observations through
$n_j$. The wavelet spectrum is defined as
\begin{equation}
  \label{eq:wav-sp_def}
{\sigma}^2_{j}=\mathbb{E}[\widehat{\sigma}^2_{j}]=\mathbb{E}[W_{j,k}^2]\quad\text{for
  all $k$}\;,
\end{equation}
where the last equality holds for $M\geq K$ since in this case
$\{W_{j,k},k\in\mathbb{Z}\}$ is weakly stationary.  We obtain the following
result which provides asymptotics of the scalogram involving a finite number of
different scales at the same time. We only provide the result for $q_0\geq2$
since the case $q_0=1$ can be directly deduced from
\cite[Theorem~2]{roueff-taqqu-2009b}.
\begin{theorem}\label{thm:multiscale}
  Suppose that Assumptions~\textbf{A}(\ref{item:A1})(\ref{item:A2}) hold with
  $q_0\geq2$. Set
  $\gamma_j=2^j$ and let
  $\{(g_j)_{j\geq0},\,g_\infty\}$ be a sequence of univariate filters
  satisfying~(W-\ref{ass:w-a})--(W-\ref{ass:w-c}) with $m=1$ and $M\geq
  \delta(q_0)+K$. Then, as $j\to\infty$,
  \begin{equation}
    \label{eq:wav-sp-asymp}
    {\sigma}^2_{j} \sim q_0!  \, (f^*(0))^{q_0} \,
L_{q_0}(\widehat{g}_{\infty})\; 2^{2 j (\delta(q_0)+K)}\;.
  \end{equation}
  Let now $j=j(N)$ be an increasing sequence such that $j\to\infty$ and
  $N2^{-j}\to\infty$. Define $n_j$, $\widehat{\sigma}^2_{j}$ and
  ${\sigma}^2_{j}$ as in~(\ref{eq:nj_def}),~(\ref{eq:scalo_def})
  and~(\ref{eq:wav-sp_def}), respectively.  Then, as $N\to\infty$,
\begin{equation}
    \label{eq:scalo-multiscale-asymp}
   \left\{ n_j^{1-2d}
     \left(\frac{\widehat{\sigma}^2_{j-u}}{{\sigma}^2_{j-u}}-1\right)\right\}_{u\geq0}
\overset{\mathrm{fidi}}{\longrightarrow}
\left\{2^{(2d-1)u}\,
\frac{L_{q_0-1}(\widehat{g}_{\infty})}{q_0!\,L_{q_0}(\widehat{g}_{\infty})}
\,Z_d(1)\right\}_{u\geq0} \;.
  \end{equation}
\end{theorem}
This theorem is proved in Section~\ref{s:proofs-main}.
Note that the constants $L_{q_0}(\widehat{g}_{\infty})$ and
$L_{q_0-1}(\widehat{g}_{\infty})$ appearing in~(\ref{eq:wav-sp-asymp})
and~(\ref{eq:scalo-multiscale-asymp}) are defined by~(\ref{e:defKp}).
Here $\overset{\mathrm{fidi}}{\longrightarrow}$ means the convergence of
finite-dimensional distributions, and since the limit depends on $u$ only
through a deterministic multiplicative constant,
we obtain, as in the multivariate case, that the multiscale limit has
cross-correlations equal to 1.

As in the multivariate case, conveniently normalized, the centered multiscale
scalogram is asymptotically a fully correlated Rosenblatt process. We recover
the results of \cite{bardet:tudor:2010} where $Y$ is the Rosenblatt process
itself. In other words Theorem 4 in ~\cite{bardet:tudor:2010} roughly
corresponds here to the case $q_0=2$. The results in \cite{CTV10} correspond to
the single scale limit for any $q_0\geq2$, which indicate a limit of the
scalogram (which corresponds to a wavelet large scale analysis) similar to that
of the variogram (which corresponds to a small scale analysis using discrete
variations).

\section{Estimation of the long memory parameter}
\label{sec:estim-long-memory}
We now consider the estimation of the long memory parameter of the observed
process $\{Y_t\}_{t\in\mathbb{Z}}$ under the assumptions of
Theorem~\ref{thm:multiscale}, that are supposed to hold all along this
section.  As already mentioned, $\{H_{q_0}(X_t)\}_{t\in\mathbb{Z}}$ has long
memory parameter $\delta(q_0)$. By~(\ref{e:K}), applying the setting of
\cite{moulines:roueff:taqqu:2007:jtsa} for dealing with processes with
stationary $K$-th increments, we get that $\{Y_t\}_{t\in\mathbb{Z}}$ itself has
long memory parameter
\begin{equation}
  \label{eq:d0-Y}
d_0=\delta(q_0)+K \;.
\end{equation}
We want to estimate this parameter from a sample $Y_1,\dots,Y_N$. A typical
wavelet estimator of $d_0$ reads
\begin{equation}
  \label{eq:d-est}
\widehat{d_0}=\sum_{i=0}^{p-1} w_{i} \log\widehat{\sigma}^2_{j+i} \;,
\end{equation}
where $w_0,\dots,w_{p-1}$ are weights such that $w_0+\dots+w_{p-1}= 0$, and
$\sum_{i=0}^{p-1} i\, w_i=1/(2\log2)$, see \cite{roueff-taqqu-2009b}. Indeed,
for this choice of weights and using~(\ref{eq:wav-sp-asymp})
and~(\ref{eq:d0-Y}), we see that, as $j\to\infty$,
  \begin{align}
\nonumber
\sum_{i=0}^{p-1} w_{i} \log{\sigma}^2_{j+i} & =
\sum_{i=0}^{p-1} w_{i} \log\left(q_0!  \, (f^*(0))^{q_0} \,
L_{q_0}(\widehat{g}_{\infty})\; 2^{2 j d_0} \right)
+ d_0 \left(\sum_{i=0}^{p-1} i\, w_{i}\right) 2 \log 2 + o(1) \\
    \label{eq:est-approx}
&= d_0 + o(1) \;.
\end{align}
Replacing ${\sigma}^2_{j+i}$ by  $\widehat{\sigma}^2_{j+i}$ in the left-hand
side of this approximation, we thus obtain an estimator $\widehat{d_0}$ of $d_0$.

To obtain the asymptotic behavior of $\widehat{d_0}$ as $j$ and $N$ go to
infinity, we first evaluate the bias, which is related to the approximation
error in the equivalence~(\ref{eq:wav-sp-asymp}). To this end, we must specify
the convergence of $f^\ast(\lambda)$ to $f^\ast(0)$ as $\lambda\to0$. A
standard assumption in the semi-parametric setup is
$$
|f^*(\lambda)-f^*(0)|\leq C f^*(0)\,|\lambda|^\beta\,\quad \lambda\in(-\pi,\pi)\;,
$$
where $\beta$ is some smoothness exponent in $(0,2]$. However here $f$ is the
spectral density of the original Gaussian process $\{X_t\}$, hence we cannot
apply directly the bound
$$
|{\sigma}^2_{j}-C_12^{2dj}|\leq C_2\,2^{(2d-\beta)j}\;,
$$
which corresponds to Relation~(26) in Theorem~1 of
\cite{moulines:roueff:taqqu:2007:jtsa} (with different notation for constants
$C_1$ and $C_2$). In fact such a bound would contradict~(\ref{eq:wav-sp-asymp})
since $d\neq d_0$, see~(\ref{e:ldparamq}) and~(\ref{eq:d0-Y}).  We must instead
work with the (generalized) spectral density, say $\tilde{f}$, of the observed
process $\{Y_t\}$.  Applying Lemma~4.1 in
\cite{clausel-roueff-taqqu-tudor-2011a}, we have that the generalized
spectral density $\tilde{f}$ of the process
$\{Y_t\}=\{\Delta^{-K}H_{q_0}(X_t)\}$ satisfies
$$
\tilde{f}(\lambda)=q_0!\;|1-\rme^{-\rmi \lambda}|^{-2 K}\,
\underset{\text{$q_0$ times}}{\underbrace{f\star\dots\star f}}(\lambda)\;,
$$
where $\star$ denotes periodic convolution. Now, by Lemma~8.2 in
\cite{clausel-roueff-taqqu-tudor-2011a}, we get
$$
q_0!\;f\star\dots\star f(\lambda) =
|1-\rme^{-\rmi \lambda}|^{-2 \delta(q_0)}
\tilde{f}^*(\lambda)\;,
$$
where $\tilde{f}^*$ denotes a nonnegative periodic function, continuous and
positive at the origin, such that
$$
|\tilde{f}^*(\lambda)-\tilde{f}^*(0)|\leq C
\tilde{f}^*(0)\,|\lambda|^{\tilde{\beta}}\,\quad \lambda\in(-\pi,\pi)\;,
$$
where $\tilde{\beta}$ is any positive number such that
$\tilde{\beta}<2\delta(q_0)$ and $\tilde{\beta}\leq\beta$.
Hence, we finally obtain
$$
\tilde{f}(\lambda)=|1-\rme^{-\rmi \lambda}|^{-2 d_0} \tilde{f}^*(\lambda) \;,
$$
and we may now apply Theorem~1 of
\cite{moulines:roueff:taqqu:2007:jtsa} and use~(\ref{eq:wav-sp-asymp}),
to obtain that
$$
\left|{\sigma}^2_{j} -q_0! (f^*(0))^{q_0} \,
L_{q_0}(\widehat{g}_{\infty})\; 2^{2 j d_0}\right|
\leq C'\,2^{j (2d_0-\tilde{\beta})}\;.
  $$
This yields
  \begin{equation}
    \label{eq:wav-sp-approx}
\left|\sum_{i=0}^{p-1} w_{i} \log{\sigma}^2_{j+i} - d_0 \right|
=O\left(2^{-\tilde{\beta}j}\right)\;,
\end{equation}
which is a more precise approximation than~(\ref{eq:est-approx}).
Observe now that $\widehat{d_0}$ given in~(\ref{eq:d-est})
satisfies the identity
$$
\widehat{d_0}
= d_0 +
\sum_{i=0}^{p-1} w_{i} \log
\left\{ 1+ \left(\frac {\widehat{\sigma}^2_{j+i}}{{\sigma}^2_{j+i}}-1\right)\right\}
+\sum_{i=0}^{p-1} w_{i} \log{\sigma}^2_{j+i} - d_0\;.
$$
Expanding $\log(x)$ in the neighborhood of $x=1$ and
using~(\ref{eq:scalo-multiscale-asymp}) and~(\ref{eq:wav-sp-approx}), we obtain
the following result.
\begin{theorem}\label{thm:estd}
  Suppose that the assumptions of Theorem~\ref{thm:multiscale} hold. As
  $N\to\infty$, if $j=j(N)$ is such
  that $j\to\infty$ and $N2^{-j}\to\infty$, then
$$
\widehat{d_0}
= d_0 + n_j^{2d-1} O_P(1)+O\left(2^{-\tilde{\beta}j}\right)\;.
$$
Moreover the $O_P$-term converges in distribution to the Rosenblatt
variable
\begin{equation}
  \label{eq:limit-est-def}
\left(\sum_{i=0}^{p-1} w_{i}2^{(1-2d)i}\right)\,
\frac{L_{q_0-1}(\widehat{g}_{\infty})}{q_0!\,L_{q_0}(\widehat{g}_{\infty})}
\,Z_d(1) \;.
\end{equation}
\end{theorem}
To optimize the asymptotic term~(\ref{eq:limit-est-def}), one should choose
weights $w_0,\dots,w_{p-1}$ which minimize the constant in parentheses. It is
interesting to note that this constant vanishes for some well-chosen weights,
but that such a choice depends on the (unknown) parameter $d$.  Observe also
that the constant approaches 0 as $d$ approaches $1/2$, since $\sum_iw_i=0$.

\begin{remark}
  To our knowledge, the non-linear semiparametric setting has not been
  considered before in this context. The closest reference appears to be
  \cite{giraitis:taqqu:1999}, where the parametric Whittle estimator is studied
  for non-linear subordinated Gaussian processes.  The comparison is difficult
  since, in the parametric approach of \cite{giraitis:taqqu:1999}, the
  asymptotic results depend on the parameterization of the spectral density
  (essentially through the two constants $\rho_1$ and $\rho_2$ defined in
  \cite{giraitis:taqqu:1999}). However similarities can be observed in these
  results: the limit can be Rosenblatt, in which case the usual $n^{-1/2}$
  parametric rate of convergence is replaced by $n^{2d-1}$, see
  \cite[Theorem~3.1]{giraitis:taqqu:1999} in the case $\rho_1=0$ and
  $\rho_2\neq0$. This situation can be compared to Theorem~\ref{thm:estd}
  above, where the limit is also Rosenblatt and the usual $n_j^{-1/2}$
  semiparametric rate is replaced by $n_j^{2d-1}$. We thus expect that a
  semiparametric Whittle approach would have an asymptotic behavior similar to
  that of $\widehat{d_0}$.
\end{remark}
% A similar conclusion holds for the wavelet Whittle estimator whose asymptotic
% behavior can be derived from that of the scalogram by following the same lines
% as in the linear case in \cite[Theorem~5]{roueff-taqqu-2009b}.

\section{Chaos expansion of the scalogram}\label{s:expansion}
Here we take $m=1$ without loss of generality, since the case $m\geq2$ can be
deduced by applying the case $m=1$ to each entry.  The purpose of this section
is to consider the scalogram $S_{n,j}$ defined in~(\ref{e:defsnj}).  and
express it as a sum of multiple integrals $\widehat{I}(\cdot)$ (defined in
Appendix A) with respect to the Gaussian random measure $\widehat{W}$. Our main
tool will be the product formula for multiple Wiener-It\^o integrals. In view
of~(\ref{eq:Wjk_representation}), $W_{j,k}$ is a multiple integral of order
$q_0$ of some kernel $f_{j,k} $, that is
\begin{equation}\label{e:a2}
W_{j,k} =\widehat{I}_{q_0} (f_{j,k} ).
\end{equation}
Now, using the product formula for multiple stochastic
integrals~(\ref{EqProdIntSto}), one gets, as shown in
Proposition~\ref{pro:decompSnjLD} that, for any
$(n,j)\in\mathbb{N}^2$,
\begin{equation}\label{e:Snjdec}
S_{n,j}-\mathbb{E}(S_{n,j})
=\frac{1}{n}\sum_{k=0}^{n-1}W_{j,k}^2-\mathbb{E}[W_{j,0}^2]=
\sum_{p=0}^{q_0-1}p! {{q_0}\choose{p}}^2\; S_{n,j}^{(p)}
\end{equation}
where, for all $0\leq p\leq q_0-1$,
\begin{equation*}
S_{n,j}^{(p)}= \widehat{I}_{2q_0-2p}(g_p)\;.
\end{equation*}
That is, for every $j,n$, the random variable  $S_{n,j}^{(p)}$ is an element of the chaos of order $2q_0-2p$.  The function $g_p(\xi)$,
$\xi=(\xi_1,\dots,\xi_{2q_0-2p})\in\mathbb{R}^{2q_0-2p}$ is
defined for every $p\in\{0,\cdots,q_0-1\}$ as
\begin{equation}\label{e:g-def}
g_p(\xi)= \frac{1}{n}\sum_{k=0}^{n-1} \left(
f_{j,k}\overline{\otimes}_{p} f_{j,k}\right)\;,
\end{equation}
where the contraction  $\overline{\otimes}_{p}$ is defined
in~(\ref{e:times-p}).

Let us formalize the above decomposition of $S_{n,j}$ and  give a
more explicit expression for the function $g_p$ in~(\ref{e:g-def}).
\begin{proposition}\label{pro:decompSnjLD}
For all non--negative integer $j$, $\{W_{j,k}\}_{k\in\mathbb{Z}}$ is a weakly stationary sequence. Moreover,
for any $(n,j)\in\mathbb{N}^2$,
\begin{equation}\label{e:decompSnjLD} S_{n,j}-\mathbb{E}(S_{n,j}) = \sum_{p=0}^{q_0-1}p!
{{q}\choose{p}}^2\; S_{n,j}^{(p)}\;,
\end{equation}
where, for all $0\leq p\leq q_0-1$,
\begin{equation}\label{e:Snjrrp1}
S_{n,j}^{(p)}= \widehat{I}_{2q_0-2p}(g_p)\;,
\end{equation}
and where, for all
  $\xi=(\xi_1,\dots,\xi_{2q_0-2p})\in\mathbb{R}^{2q_0-2p}$,
\begin{multline}\label{e:bg}
g_p(\xi)= D_n(\gamma_j\left(\xi_1+\dots+\xi_{2q_0-2p}\right))\\
\times
\prod_{i=1}^{2q_0-2p}[\sqrt{f(\xi_i)}\1_{(-\pi,\pi)}(\xi_i)]
\times
\widehat{\kappa}_j^{(p)}(\xi_1+\dots+\xi_{q_0-p},\xi_{q_0-p+1}+\dots+\xi_{2q_0-2p})
\;.
\end{multline}
Here $f$ denotes the spectral density (\ref{e:sdf}) of the underlying Gaussian
process $X$ and
\begin{equation}\label{eq:dirichlet}
D_n(u)=\frac{1}{n}\sum_{k=0}^{n-1}\rme^{\rmi k u}=
\frac{1-\rme^{\rmi n u}}{n(1-\rme^{\rmi  u})}\;,
\end{equation}
denotes the normalized Dirichlet kernel. % Note that the last equality holds for $u\notin 2\pi \mathbb{Z}$ and for every $u\in 2\pi \mathbb{Z}$ we have $D_{n}(u)=1$.
Finally, for
$\xi=(\xi_1,\xi_2)\in\mathbb{R}^2$, if $p\neq 0$,
\begin{equation}\label{e:intrepKjp1}
\widehat{\kappa}_{j}^{(p)}(\xi_1,\xi_2)=\int_{(-\pi,\pi)^{p}}\left(\prod_{i=1}^{p}f(\lambda_i)\right)\;
\widehat{h}_j^{(K)}(\lambda_1+\dots+\lambda_p+\xi_1)\overline{\widehat{h}_{j}^{(K)}(\lambda_1+\dots+\lambda_p-\xi_2)}\;\rmd^p\lambda\;,
\end{equation}
and, if $p= 0$,
\begin{equation}\label{e:intrepKjp1bis}
\widehat{\kappa}_{j}^{(p)}(\xi_1,\xi_2)=\widehat{h}_j^{(K)}(\xi_1)\overline{\widehat{h}_j^{(K)}(\xi_2)}.
\end{equation}
\end{proposition}

\noindent{\bf Notation.}
In~(\ref{e:intrepKjp1}), $\rmd^p\lambda$ refers to $p$-dimensional Lebesgue
measure integration.  To simplify the notation, we shall denote by $\Sigma_q$,
the $\mathbb{C}^{q}\to\mathbb{C}$ function defined, for all $q\in\mathbb{Z}_+$
and $y=(y_1,\dots,y_{q})\in\mathbb{C}^{q}$, by
\begin{equation}\label{e:DefPartialSums1}
\Sigma_q(y)=\sum_{i=1}^{q} y_i\;,
  \end{equation}
and for any $(q_1,q_2)\in\mathbb{Z}_+^2$, we denote by
$\Sigma_{q_1,q_2}$  the $\mathbb{C}^{q_1}\times
\mathbb{C}^{q_2}\to\mathbb{C}^2$ function defined for all
$y=(y_1,\dots,y_{q_1+q_2})\in\mathbb{C}^{q_1}\times
\mathbb{C}^{q_2}$ by
\begin{equation}\label{e:DefPartialSums2}
\Sigma_{q_1,q_2}(y)=\left(\sum_{i=1}^{q_1} y_i,\sum_{i=q_1+1}^{q_2}y_i\right)\;.
  \end{equation}

  With these notations, (\ref{e:Snjrrp1}),~(\ref{e:intrepKjp1})
  and~(\ref{e:intrepKjp1bis}) become respectively
\begin{align}\label{e:Snjrrp}
&S_{n,j}^{(p)}=\widehat{I}_{2q_0-2p}\left(
D_n\circ\Sigma_{2q_0-2p}(\gamma_j\times\cdot)\times
[\sqrt{f}\1_{(-\pi,\pi)}]^{\otimes(2q_0-2p)}
\times \widehat{\kappa}_j^{(p)}\circ\Sigma_{q_0-p,q_0-p}\right)\;,\\
\label{e:intrepKjp}
&\widehat{\kappa}_{j}^{(p)}(\xi_1,\xi_2)=
\begin{cases}
\int_{(-\pi,\pi)^{p}}
f^{\otimes p}(\lambda)\;
\widehat{h}_j^{(K)}(\Sigma_p(\lambda)+\xi_1)\overline{\widehat{h}_{j}^{(K)}(\Sigma_p(\lambda)-\xi_2)}\;\rmd^p\lambda&\text{if
}p\neq 0,\\
[\widehat{h}_j^{(K)}\otimes\overline{\widehat{h}_j^{(K)}}](\xi_1,\xi_2)
&\text{if
}p=0\;,
\end{cases}
\end{align}
where $\circ$ denotes the composition of functions,
$\lambda=(\lambda_1,\cdots,\lambda_p)$ and $f^{\otimes
  p}(\lambda)=f(\lambda_1)\cdots f(\lambda_p)$ is written as a tensor
product.

\begin{remark}
  The kernel $\widehat{\kappa}_{j}^{(p)}$ can also be expressed in terms of the
  the covariance sequence of the process $X$, namely,
\begin{equation}\label{e:exphjkp}
\widehat{\kappa}_{j}^{(p)}(\xi_1,\xi_2)=\sum_{m\in\mathbb{Z}^2}
h_{j}^{(K)}(m_1)h_{j}^{(K)}(m_2)\;\mathbb{E}(X_{m_2}X_{m_1})^p \;\rme^{-\rmi
(m_1\xi_1+ m_2\xi_2)}\;.
\end{equation}
This follows from the relation
\begin{equation*}
\mathbb{E}(X_{m_{2}}X_{m_{1}})=\int_{-\pi }^{\pi}
e^{i(m_{2}-m_{1})\lambda } f(\lambda) \rmd\lambda\;,
\end{equation*}
and  (\ref{e:HjkVSHj}) and the definition~(\ref{e:dF}) of
the discrete Fourier transform $\widehat{h}_{j}$.
\end{remark}

\noindent{\bf Proof of Proposition~\ref{pro:decompSnjLD}}
By~(\ref{e:defsnj}),
\begin{eqnarray}\label{e:a3}
S_{n,j}&=&\frac{1}{n}\sum_{k=0}^{n-1}W_{j,k}^2\;.
\end{eqnarray}
Using ~(\ref{e:a2}) and the product formula for multiple stochastic integrals ~(\ref{EqProdIntSto})
of Proposition~\ref{pro:intproduct}, we have
\begin{equation}\label{e:prodwjk}
W_{j,k}^2=\widehat{I}_{q_0}(f_{j,k})\widehat{I}_{q_0}(f_{j,k})
=\sum_{p=0}^{q_0} p!
{{q_0}\choose{p}}^2
\widehat{I}_{2q_0-2p}\left(f_{j,k} \overline{\otimes}_p
f_{j,k}\right)\;.
\end{equation}
Therefore,
\begin{equation}\label{e:a4}
S_{n,j}=\frac{1}{n}\sum_{k=0}^{n-1}W_{j,k}^2=\sum_{p=0}^{q_0}
p!{{q_0}\choose{p}}^2\widehat{I}_{2q_0-2p}\left(g_p\right)\;,
\end{equation}
where
\[
g_{p}=\frac{1}{n}\sum_{k=0}^{n-1}f_{j,k}\overline{\otimes}_p
f_{j,k}\;.
\]
By~(\ref{e:fjk_representation1}), for all
$\xi=(\xi_1,\cdots,\xi_{q_0})\in\mathbb{R}^{q_0}$,
\begin{equation}\label{e:fjk_representation2}
f_{j,k}(\xi)=
\exp\circ\Sigma_{q_0}(\rmi k\gamma_j\xi)\left(\widehat{h}_{j}^{(K)}\circ\Sigma_q(\xi)\right)\left(f^{\otimes q_0}(\xi)\right)^{1/2}
\1_{(-\pi,\pi)}^{\otimes q_0}(\xi)\;.
\end{equation}

If, $p=1,2, \ldots, q_0 -1$, let $\xi=(\xi_1,\cdots,\xi_{2q_0-2p})$. The
contraction $f_{j,k} \overline{\otimes}_p f_{j,k}$ defined on
$\mathbb{R}^{2q_0-2p}$ equals by~(\ref{e:times-p}),
\begin{eqnarray*}
&&f_{j,k}
\overline{\otimes}_p f_{j,k}(\xi)\\
&=&\int_{\mathbb{R}^p}f_{j,k}(\xi_1,\cdots,\xi_{q_0-p},s)f_{j,k}(\xi_{q_0-p+1},\cdots,\xi_{2q_0-2p},-s)\rmd^p s\\
&=&
\exp\circ\Sigma_{2q_0-2p}(\rmi k\gamma_j\xi)
\times [\sqrt{f}\1_{(-\pi,\pi)}]^{\otimes 2q_0-2p}(\xi)\\
&&\times
\int_{\mathbb{R}^p}
\widehat{h}_{j}^{(K)}(\xi_1+\dots+\xi_{q_0-p}+\Sigma_p(\lambda))
\widehat{h}_{j}^{(K)}(\xi_{q_0-p+1}+\dots+\xi_{2q_0-2p}-\Sigma_p(\lambda))
\times [f\1_{(-\pi,\pi)}]^{p}(\lambda)
\;\rmd^p\lambda\\
&=& \exp\circ\Sigma_{2q_0-2p}(\rmi k\gamma_j\xi) \times
[\sqrt{f}\1_{(-\pi,\pi)}]^{\otimes 2q_0-2p}(\xi) \times
\widehat{\kappa}_j^{(p)}\circ\Sigma_{q_0-p,q_0-p}(\xi)\;,
\end{eqnarray*}
where $\widehat{\kappa}_{j}^{(p)}$ is defined
by~(\ref{e:intrepKjp1}),(\ref{e:intrepKjp1bis}), or equivalently
by~(\ref{e:intrepKjp}),(\ref{e:intrepKjp1bis}) and where we used
that
$\overline{\widehat{h}_{j}^{(K)}(\cdot)}=\widehat{h}_{j}^{(K)}(-\cdot)$.
We therefore get that $g_p$ is a function with $2q_0-2p $ variables
given by
\begin{equation*}
g_p(\xi) = \frac{1}{n}
\sum_{k=0}^{n-1}\exp\circ\Sigma_{2q_0-2p}(\rmi k\gamma_j\xi)
\times [\sqrt{f}\1_{(-\pi,\pi)}]^{\otimes 2q_0-2p}(\xi) \times
\widehat{\kappa}_j^{(p)}\circ\Sigma_{q_0-p,q_0-p}(\xi)\;.
\end{equation*}

The Dirichlet kernel $D_{n}$ appears when one computes the sum
$\frac{1}{n} \sum_{k=0}^{n-1}\exp\circ\Sigma_{2q_0-2p}(\rmi
k\gamma_j\xi)$. This implies the  formula (\ref{e:bg}).

The chaos of order zero does not appears in (\ref{e:decompSnjLD})
where $S_{n,j}-\mathbb{E}(S_{n,j})$ is considered. It appears
however in the expression~(\ref{e:a4}) of $S_{n,j}$ in the term
with $p=q_0$ where $\widehat{I}_{2q_0-2p}=\widehat{I}_0$. In this
case, we have
\[
q_0!\widehat{I_0}(f_{j,k}\overline{\otimes}_p
f_{j,k})=
q_0!\|f_{j,k}\|_{L^{2}(\mathbb{R}^{q_0})}^2=\mathbb{E}(|W_{j,k}|^2)\;,
\]
corresponding in~(\ref{e:a4}) to the deterministic term
\[
\frac{1}{n}\sum_{k=1}^{n}\mathbb{E}(|W_{j,k}|^2)
=\mathbb{E}(|W_{j,0}|^2) =
\mathbb{E}(S_{n,j})\;,
\]
by~(\ref{e:a3}). Therefore
$S_{n,j}-\mathbb{E}(S_{n,j})$ can be expressed
as~(\ref{e:decompSnjLD}). \null\hfill $\Box$\par\medskip

As we can see from~(\ref{e:decompSnjLD}), the random variable
$S_{n,j}$ can be expanded into a sum of multiple stochastic
integrals starting from order zero (which corresponds to the
deterministic term $\mathbb{E}(S_{n,j})$). The order of the chaos
appearing in the decomposition of $S_{n,j}$ could be greater or
smaller than the critical value $1/(1-2d)$. This means that
$S_{n,j}$ may admit summands with long-range dependence (orders
smaller than $1/(1-2d)$ ) and short-range dependence (orders
greater than $1/(1-2d)$). We will see that these two kind of terms
have different behavior.  Another issue concerns $p$, the order of
the contraction in the product formula for multiple integrals. The
case $p=0$ must be discussed separately because the function
$\widehat{\kappa}^{(p)}_{j}$ in~(\ref{e:intrepKjp1}) has the
special form~(\ref{e:intrepKjp1bis}) if $p=0$.

To study $S_{n,j}$ as $j,n\to\infty$, we need to study
$S_{n,j}^{(p)}$ which is given in~(\ref{e:Snjrrp}). We first
estimate the $L^2$ norm of $S_{n,j}^{(p)}$.

\section{An upper bound for the $L^2$ norm of the terms $S_{n,j}^{(p)}$}\label{s:L2bound}
To identify the leading term of the sum $S_{n,j}-\mathbb{E}(S_{n,j})$, we will
give an upper bound for the $L^2$ norms of the terms $S_{n,j}^{(p)}$ $0\leq
p<q_0$ defined in~(\ref{e:Snjrrp1}) and~(\ref{e:Snjrrp}). Then, in
Section~\ref{s:leadingterms}, we investigate the asymptotic behavior of the
leading term of $S_{n,j}$. It directly implies the required result about the
asymptotic bahavior of the scalogram.  The expression~(\ref{e:Snjrrp}) of
$S_{n,j}^{(p)}$ involves the kernel $\widehat{\kappa}^{(p)}_j$
in~(\ref{e:intrepKjp}) which vanishes when $\xi_1=0$ or $\xi_2=0$ if $p=0$
because $\widehat{h}_j(0)=0$ by~(\ref{e:majoHj}). But the
expression~(\ref{e:intrepKjp}) of $\widehat{\kappa}_j^{(p)}$ implies that it
does not vanish if $p>0$ because
\[
\widehat{\kappa}_j^{(p)}(0,0)=\int_{(-\pi,\pi)^p}\left(\prod_{i=1}^p
f(\lambda_i)\right)\left|\widehat{h}_j(\Sigma_p(\lambda))\right|^2\rmd^p\lambda>0\;.
\]
All these considerations lead
one to distinguish the following two cases :
\begin{itemize}
\item The case $p\neq 0$.
\item The case $p= 0$.
\end{itemize}
As for the Rosenblatt process considered by~\cite{bardet:tudor:2010},
the case $p=0$ requires different bounds and thus must be treated separately.

\subsection{The case $p\neq 0$}

Recall the expansion~(\ref{e:decompSnjLD}). In the case $p\neq 0$ we now give an
upper bound of $\|S_{n,j}^{(p)}\|_{2}=\mathbb{E}(|S_{n,j}^{(p)}|^2)^{1/2}$ with
$0<p< q_0<1/(1-2d)$.

\begin{proposition}\label{pro:pneqzSD}Let $0<p<q_0<1/(1-2d)$.
There exists some $C>0$ whose value depends only on $p,d,q_0$ and $f^*$
such that for all $n,j\geq2$
\begin{equation}\label{e:boundL2P51}
\|S_{n,j}^{(p)}\|_{2}\leq C (\log n)^\varepsilon
\;n^{-\min(1-2\delta(q_0-p),1/2)}\,\gamma_j^{2\delta(q_0)+2K}\;,
\end{equation}
where $\varepsilon=1$ if $\delta(q_0-p)=1/4$ and  $\varepsilon=0$ otherwise.
\end{proposition}
\begin{proof}
Let $C,C_1,\cdots$ be positive constants that may change from line
to line. Set $r=q_0-p \geq 1$. We perform the change of variable $y=n\gamma_j\xi$ in the integral
  expression of $S_{n,j}^{(p)}$ given by~(\ref{e:Snjrrp}) and deduce that
\begin{eqnarray*}
\mathbb{E}\left|S_{n,j}^{(p)}\right|^2&=&
\frac{1}{(n\gamma_j)^{2r}}\int_{\mathbb{R}^{2r}}
\left|D_n\circ\Sigma_{2r}\left(\frac{y}{n}\right)\right|^2
\;\left(\prod_{i=1}^{2r}(f\1_{(-\pi,\pi)})(\frac{y_i}{n\gamma_j})\right)
\;\left|\widehat{\kappa}_{j}^{(p)}\circ\Sigma_{r,r}\left(\frac{y}{n\gamma_j}\right)\right|^2
\;\rmd^{2r} y\;.
\end{eqnarray*}
We now use the expression of $f$ given by~(\ref{e:sdf}), the
boundedness of $f^*$, the bound of Dirichlet kernel given by
Lemma~\ref{lem:Dn} and the bound of $\widehat{\kappa}_{j}^{(p)}$
given by Lemma~\ref{lem:suphjkp}. Hence one deduces that there
exists some $C_1>1$ depending only on $p,d$ such that
\begin{equation}\label{e:boundSnj51}
\mathbb{E}\left|S_{n,j}^{(p)}\right|^2\leq C_1\gamma_j^{-2r(1-2d)}\gamma_j^{4(\delta(p)+K)}I_{n,j}=C_1\gamma_j^{-2+4\delta(r)+4\delta(p)}\gamma_j^{4K}I_{n,j}\;,
\end{equation}
where
\[
I_{n,j}=\int_{(-n\gamma_j\pi,n\gamma_j\pi)^{2r}}
\frac{n^{-2r(1-2d)}\left|g\circ
\Sigma_{r,r}(\frac{y}{n\gamma_j})\right|^2\;\rmd^{2r} y}
{\left(1+n\left|\{\Sigma_{2r}(n^{-1}y)\}\right|\right)^2\prod_{i=1}^{2r}\left|y_i\right|^{2d}}\;,
\]
with
\[
g(z_1,z_2)=\frac{1}{(1+\gamma_j|\{z_1\}|)^{\delta(p)}(1+\gamma_j|\{z_2\}|)^{\delta(p)}}\;.
\]
We now bound the integral $I_{n,j}$. To this end, perform the
successive change of variables
\[
u_1=\frac{y_1+\cdots+y_{r}}{n},\,\cdots,\,u_{r}=\frac{y_{r}}{n},v_1=\frac{y_{r+1}+\cdots+y_{2r}}{n},\,\cdots,\,v_{r}=\frac{y_{2r}}{n}\;,
\]
so that
\begin{eqnarray*}
y_i=n(u_i-u_{i+1})\mbox{ for }1\leq i\leq r-1,\; y_r=n u_r,\\
y_i=n(v_{i-r}-v_{i-r+1})\mbox{ for }r+1\leq i\leq 2r-1,\;
y_{2r}=n u_{r}\;.
\end{eqnarray*}
In addition, observe that for any
$m\in\mathbb{Z}_+\setminus\{0\}$,
$\left(y_1,\cdots,y_m\right)\in(-n\gamma_j\pi,n\gamma_j\pi)^m\;,$
implies that
$y_1+\cdots+y_m\in(-m(n\gamma_j)\pi,m(n\gamma_j)\pi)$. Hence,
there exists some constant $C$ depending only on $r,d$ such
that
\begin{equation}
  \label{eq:ineg62}
I_{n,j}\leq C\int_{-\gamma_j\pi r}^{\gamma_j\pi
r}\int_{-\gamma_j\pi r}^{\gamma_j\pi
r}\frac{J_{r,\gamma_j\pi}(u_1;2d\mathrm{1}_r)
J_{r,\gamma_j\pi}(v_1;2d\mathrm{1}_r)\rmd u_1\rmd
v_1} {(1+n\left|\{u_1+v_1\}\right|)^{2}
(1+\gamma_j\left|\{\frac{u_1}{\gamma_j}\}\right|)^{2\delta(p)}(1+\gamma_j\left|\{\frac{v_1}{\gamma_j}\}\right|)^{2\delta(p)}}\;,
\end{equation}
where we used the definition of $J_{m,a}(s;\beta)$ in Lemma~\ref{lem:Ja} with
the notation $\mathrm{1}_r$ for the $r$-dimensional vector with all entries
equal to 1, that is, we set $m=r$, $a=\gamma_j\pi$, $\beta_1=\dots=\beta_m=2d$
in~(\ref{e:Jma}).  We now apply Lemma~\ref{lem:Ja}. Since
$m=r<1/(1-2d)$, we are in Case~(i) and we get that there exists some $C>0$
depending only on $r,d$ such that
\[
J_{r,\gamma_j\pi}(s;2d\mathrm{1}_r)
\leq C|s|^{-2\delta(r)}\,\;\mbox{
for all }s\in\mathbb{R}\;.
\]
Then there exists some constant $C_2>1$ depending only on $r,d$ such
that
\begin{equation}
  \label{eq:ineg63}
I_{n,j}\leq C_2
\int_{-\gamma_j\pi r}^{\gamma_j\pi r}\int_{-\gamma_j\pi r}^{\gamma_j\pi r}\frac{|u_1|^{-2\delta(r)}|v_1|^{-2\delta(r)}\rmd
u_1\;\rmd v_1}
{(1+n\left|\{u_1+v_1\}\right|)^2\left(1+\gamma_j\left|\{\frac{u_1}{\gamma_j}\}\right|\right)^{2\delta(p)}\;
\left(1+\gamma_j\left|\{\frac{v_1}{\gamma_j}\}\right|\right)^{2\delta(p)}}\;.
\end{equation}
Now use the inequality $|\{x\}|\leq |x|$ valid on $x\in\mathbb{R}$.
Since $\delta(r)\geq 0$,
\[
I_{n,j}\leq
C_2
\int_{-\gamma_j\pi r}^{\gamma_j\pi r}\int_{-\gamma_j\pi r}^{\gamma_j\pi r}\frac{|\gamma_j\{\frac{u_1}{\gamma_j}\}|^{-2\delta(r)}|\gamma_j\{\frac{v_1}{\gamma_j}\}|^{-2\delta(r)}\rmd
u_1\rmd v_1} {(1+n\left|\{u_1+v_1\}\right|)^2
\left(1+\gamma_j\left|\{\frac{u_1}{\gamma_j}\}\right|\right)^{2\delta(p)}\;
\left(1+\gamma_j\left|\{\frac{v_1}{\gamma_j}\}\right|\right)^{2\delta(p)}}\;.
\]
By $2\pi$--periodicity of $x\mapsto \{x\}$, the integrand is
$(2\gamma_j\pi)$-periodic with respect to both variables $u_1$ and $v_1$ and we
get that
\begin{equation}
  \label{eq:ineg64}
I_{n,j}\leq
C_3
\int_{-\gamma_j\pi }^{\gamma_j\pi }\int_{-\gamma_j\pi
}^{\gamma_j\pi
}\frac{|u_1|^{-2\delta(r)}|v_1|^{-2\delta(r)}\rmd u_1\;\rmd
v_1} {(1+n\left|\{u_1+v_1\}\right|)^2(1+|u_1|)^{2\delta(p)}
(1+|v_1|)^{2\delta(p)}}\;.
\end{equation}
To deal with the fractional parts, we now partition
$(-\gamma_j\pi,\gamma_j\pi)^2$ using the following domains
\[
\Delta_{j}^{(s)}=\{(u_1,v_1)\in
(-\gamma_j\pi,\gamma_j\pi)^2,|u_1+v_1-2\pi s|\leq \pi\}\;,
\]
with $s\in\{-\gamma_j,\dots,\gamma_j\}$, so that $I_{n,j}=A+2B$ with
\[
A=\int_{\Delta_{j}^{(0)}
}\frac{|u_1|^{-2\delta(r)}|v_1|^{-2\delta(r)}\rmd u_1\;\rmd
v_1} {(1+n\left|u_1+v_1\right|)^2(1+|u_1|)^{2\delta(p)}
(1+|v_1|)^{2\delta(p)}}\;,
\]
and
\[
B=\sum_{s=1}^{\gamma_j}\int_{\Delta_{j}^{(s)}
}\frac{|u_1|^{-2\delta(r)}|v_1|^{-2\delta(r)}\rmd u_1\;\rmd
v_1} {(1+n\left|u_1+v_1-2\pi s\right|)^2(1+|u_1|)^{2\delta(p)}
(1+|v_1|)^{2\delta(p)}}\;.
\]
Let us now bound separately $A$ and $B$.  To bound $A$, we distinguish two
cases : $4\delta(r)>1$ and $4\delta(r)\leq 1$. In the first case, observe that
$(1+|u|)^{2\delta(p)}\geq 1$ holds on $\mathbb{R}$ and perform the change of
variables $u'_1=n u_1$ and $v'_1=n v_1$. Then
\begin{equation}\label{e:boundIP51a}
A\leq
n^{-2+4\delta(r)}\int_{\mathbb{R}^2}\frac{|u'_1|^{-2\delta(r)}|v'_1|^{-2\delta(r)}\rmd
u'_1\;\rmd v'_1} {(1+\left|u'_1+v'_1\right|)^2}\leq C
n^{-2+4\delta(r)}\;,
\end{equation}
since the integral is bounded. This follows from Lemma~8.4
of~\cite{clausel-roueff-taqqu-tudor-2011a} applied with $M_1=2$, $M_2=0$, $q=2$, $a=0$,
$\beta_1=\beta_2=2\delta(r)$.\\

In the case where $4\delta(r)\leq 1$, setting $t_1= u_1+v_1$, we get that
\[
A\leq \int_{-\pi}^{-\pi}\frac{\rmd t_1}
{(1+n|t_1|)^2}
\left[\int_{-\gamma_j\pi}^{\gamma_j\pi}\frac{|t_1-v_1|^{-2\delta(r)}|v_1|^{-2\delta(r)}\rmd
v_1}{(1+|t_1-v_1|)^{2\delta(p)}(1+|v_1|)^{2\delta(p)}}\right]\;.
\]
We now split the integral in brackets into two terms
\[
\int_{|v_1|\leq
2|t_1|}\frac{|t_1-v_1|^{-2\delta(r)}|v_1|^{-2\delta(r)}\rmd
v_1}{(1+|t_1-v_1|)^{2\delta(p)}(1+|v_1|)^{2\delta(p)}}
+\int_{2|t_1|\leq|v_1|\leq\gamma_j\pi}
\frac{|t_1-v_1|^{-2\delta(r)}|v_1|^{-2\delta(r)}\rmd
v_1}{(1+|t_1-v_1|)^{2\delta(p)}(1+|v_1|)^{2\delta(p)}}
\]
Consider the first integral. Since $4\delta(r)\leq 1$,
Lemma~\ref{lem:Ja} (case (ii) or (iv)) applied with $m=2$, $a=2|t_1|,
s_1=t_1$, $\beta_1=\beta_2=2\delta(r)$ then
implies that for some $C>0$ depending on $r,d$
\begin{eqnarray*}
\int_{|v_1|\leq
2|t_1|}\frac{|t_1-v_1|^{-2\delta(r)}|v_1|^{-2\delta(r)}\rmd
v_1}{(1+|t_1-v_1|)^{2\delta(p)}(1+|v_1|)^{2\delta(p)}}&\leq&
\int_{|v_1|\leq
2|t_1|}|t_1-v_1|^{-2\delta(r)}|v_1|^{-2\delta(r)}\rmd
v_1\\
&\leq& C|t_1|^{1-4\delta(r)}\;.
\end{eqnarray*}
Now consider the second integral. Note that $|v_1|\geq 2|t_1|$
implies $|v_1-t_1|\geq |v_1|-|t_1|\geq |v_1|/2$.
We get that
\begin{eqnarray*}
\int_{2|t_1|\leq|v_1|\leq\gamma_j\pi}
\frac{|t_1-v_1|^{-2\delta(r)}|v_1|^{-2\delta(r)}\rmd
v_1}{(1+|t_1-v_1|)^{2\delta(p)}(1+|v_1|)^{2\delta(p)}}&\leq&
C\int_{2|t_1|}^{\gamma_j\pi}\frac{|v_1|^{-2\delta(r)}|v_1|^{-2\delta(r)}\rmd
v_1}{(1+|v_1|)^{2\delta(p)}(1+|v_1|)^{2\delta(p)}}\\
&\leq&
C\int_{2|t_1|}^{\gamma_j\pi}\frac{|v_1|^{-4\delta(r)}\rmd
v_1}{(1+|v_1|)^{4\delta(p)}}\\
&=&O\left((1+|\log|t_1| |)^{\varepsilon}\right)\;,
\end{eqnarray*}
where we used that $-4\delta(r)\geq-1$ with equality if and only if
$\varepsilon=1$ and that $4(\delta(r)+\delta(p))=4\delta(q_0)+2>2$.
%  and that, for any
% $\beta>0$ and $\beta'<1$,
% \begin{equation}\label{e:majoint}
% \int_{-a}^a \frac{|x|^{-\beta'}}{(1+|x|)^\beta}\rmd
% x=O(a^{\max(1-\beta'-\beta,0)}|\log(a)|^{\varepsilon(\beta'+\beta)})\mbox{ when }a\to +\infty,\;
% \end{equation}
Hence, if $4\delta(r)\leq1$
\begin{equation}\label{e:boundIP51b}
A\leq C\left(\int_{-\pi}^\pi \frac{(1+|\log|t_1| |)^{\varepsilon}\;\rmd t_1}
{(1+n|t_1|)^2}\right)
\leq Cn^{-1} \, (\log n)^\varepsilon\;.
\end{equation}
To sum up Equations~(\ref{e:boundIP51a})
and~(\ref{e:boundIP51b}), we can write
\begin{equation}\label{e:boundIP51}
A\leq
C\, (\log n)^\varepsilon\,n^{-\min(2-4\delta(r),1)}\;.
\end{equation}
To bound $B$ observe that, on $\mathbb{R}^2$, if $|u_1|\leq
|u_1+v_1|/2$ then
\[
|v_1|=|(u_1+v_1)-u_1|\geq |u_1+v_1|-|u_1|\geq |u_1+v_1|/2\;.
\]
Hence either $|u_1|\geq |u_1+v_1|/2$ or $|v_1|\geq |u_1+v_1|/2$.
Set
\[
\Delta_{j}^{(s,1)}=\{(u_1,v_1)\in \Delta_{j}^{(s)},\,|u_1|\geq
|u_1+v_1|/2\}\;,
\]
and its symmetric set
\[
\Delta_{j}^{(s,2)}=\{(u_1,v_1)\in \Delta_{j}^{(s)},\,|v_1|\geq
|u_1+v_1|/2\}\;.
\]
Then, since $\delta(r),\delta(p)> 0$, for any
$s\in\{-\gamma_j,\cdots,-1,1,\cdots,\gamma_j\}$,
\begin{eqnarray*}
B^{(s,1)}&=&\int_{\Delta_{j}^{(s,1)}
}\frac{|u_1|^{-2\delta(r)}|v_1|^{-2\delta(r)}\rmd u_1\;\rmd
v_1} {(1+n\left|\{u_1+v_1\}\right|)^2(1+|u_1|)^{2\delta(p)}
(1+|v_1|)^{2\delta(p)}}\\
&\leq& C \int_{\Delta_{j}^{(s,1)}
}\frac{|u_1|^{-2(\delta(r)+\delta(p))}|v_1|^{-2\delta(r)}\rmd
u_1\;\rmd v_1}
{(1+n\left|\{u_1+v_1\}\right|)^2(1+|v_1|)^{2\delta(p)}}\\
&\leq& C \int_{\Delta_{j}^{(s,1)}
}\frac{|u_1+v_1|^{-2(\delta(r)+\delta(p))}|v_1|^{-2\delta(r)}\rmd
u_1\;\rmd v_1}
{(1+n\left|\{u_1+v_1\}\right|)^2(1+|v_1|)^{2\delta(p)}}\;.
\end{eqnarray*}
Setting $t_1=n(u_1+v_1)$, we get that
\[
B^{(s,1)} \leq C
n^{-1+2\delta(r)+2\delta(p)}\left(\int_{t_1=2\pi ns-\pi
n}^{2\pi ns+\pi n}\frac{|t_1|^{-2\delta(r)-2\delta(p)}\rmd
t_1} {(1+\left|t_1-2\pi n
s\right|)^2}\right)\left(\int_{-\gamma_j\pi}^{\gamma_j\pi}\frac{|v_1|^{-2\delta(r)}\rmd
v_1}{(1+|v_1|)^{2\delta(p)}}\right)\;.
\]
Set $w_1=t_1-2\pi n s$. Since $s\neq 0$, we have
\[
B^{(s,1)} \leq C
n^{-1+2\delta(r)+2\delta(p)}(n(2|s|-1))^{-2\delta(r)-2\delta(p)}\left(\int_{\mathbb{R}}(1+|w_1|)^{-2}\rmd
w_1\right)\left(\int_{-\gamma_j\pi}^{\gamma_j\pi}\frac{|v_1|^{-2\delta(r)}\rmd
v_1}{(1+|v_1|)^{2\delta(p)}}\right)\;,
\]
and the same bound holds on $B^{(s,2)}$ by symmetry. Hence
\begin{equation}\label{e:majoJs1}
B=\sum_{s=1}^{\gamma_j}(B^{(s,1)}+B^{(s,2)})\leq
Cn^{-1}\left(\sum_{|s|=1}^{\gamma_j}(2|s|-1)^{-2\delta(r)-2\delta(p)}\right)\left(\int_{-\gamma_j\pi}^{\gamma_j\pi}\frac{|v_1|^{-2\delta(r)}\rmd
v_1}{(1+|v_1|)^{2\delta(p)}}\right)\;.
\end{equation}
Using $2\delta(p)+2\delta(r)=\delta(q_0)+1>1$, we
deduce from~(\ref{e:majoJs1}) that $B=O(n^{-1})$ and, with~(\ref{e:boundIP51}),
$I_{n,j}=A+B=O((\log n)^\varepsilon\,n^{-\min(2-4\delta(r),1)})$. With~(\ref{e:boundSnj51}) and
$\delta(p)+\delta(r)=\delta(q_0)+1/2$, we obtain~(\ref{e:boundL2P51}).
\end{proof}
\subsection{The case $p=0$}\label{s:casepnull}
Here the situation is different from the previous case $p\neq 0$
since the kernel $\widehat{\kappa}_j^{(p)}$ involved in the
definition of $S_{n,j}^{(p)}$ has a different expression when
$p=0$ and vanishes when $\xi_1=0$ or $\xi_2=0$. It implies that
the bound in Proposition~\ref{pro:snjpnull} involves $n^{-1/2}$
instead of $n^{-1+\delta(q_0)}$ as could be expected from the case $p>0$ in Proposition~\ref{pro:pneqzSD}. Further, an additional
assumption on the moments of the wavelet is required which is
consistent with the results proved in the Gaussian case
in~\cite{moulines:roueff:taqqu:2007:jtsa} (corresponding to $q_0=1$)
where $M$ is assumed to be greater than $K+d$.
\begin{proposition}\label{pro:snjpnull}
Assume that $M\geq \delta(q_0)+K$. Then there
exists some $C>1$ whose values depend only on $q_0,d$ such that for
any $n,j$
\begin{equation}\label{e:boundL2P52}
\|S_{n,j}^{(0)}\|_{L^{2}(\Omega)}=\mathbb{E}(|S_{n,j}^{(0)}|^2)^{1/2}\leq
C\;n^{-1/2}\gamma_j^{2\delta(q_0)+2K}\;.
\end{equation}
\end{proposition}
\begin{proof}
We denote by $C$ a positive constant that may
change at each appearance, but whose value does neither depend on
$n$ nor $j$. Since $p=0$,
$\widehat{\kappa}^{(0)}_j=\widehat{h}^{(K)\otimes2}_j$
by~(\ref{e:intrepKjp1bis}). Then, setting $y=(n\gamma_j)^{-1}\xi$
in~(\ref{e:Snjrrp}), we get
\begin{eqnarray}\label{e:L2peq0}
&&\mathbb{E}\left|S_{n,j}^{(0)}\right|^2\\
&&=\frac{1}{(n\gamma_j)^{2q_0}}\int_{\mathbb{R}^{2q_0}}
\left|D_n\circ\Sigma_{2q_0}(\frac{y}{n})\right|^2
\;(f\1_{(-\pi,\pi)})^{\otimes(2q_0)}(\frac{y}{n\gamma_j})
\;\left|\widehat{h}^{(K)\otimes2}_j\circ\Sigma_{q_0,q_0}(\frac{y}{n\gamma_j})\right|^2\rmd^{2q_0}
y\;.\nonumber
\end{eqnarray}
We now use the bound of the Dirichlet kernel given by Lemma~\ref{lem:Dn}, the
definition of $f$ given by Equation~(\ref{e:sdf}) with the boundedness of
$f^*$, the bound of $\widehat{h}_j^{(K)}$ given by Equation~(\ref{e:hK}). Then
we deduce that
\begin{equation}\label{e:boundSnj1}
\mathbb{E}[|S_{n,j}^{(0)}|^2]\leq
C\,\gamma_j^{-2q_0(1-2d)}\gamma_j^{2(2K+1)}I_{n,j}=C\gamma_j^{4(\delta(q_0)+K)}I_{n,j}\;,
\end{equation}
where $\delta(\cdot)$ is defined by~(\ref{e:ldparamq}) and where for any $j,n$
$$
I_{n,j}=n^{-2q_0(1-2d)}\int_{(-n\gamma_j\pi,n\gamma_j\pi)^{2q_0}}
g\circ\Sigma_{q_0,q_0}(\frac{y}{n})\;\left(\prod_{i=1}^{2q_0}|y_i|^{-2d}\right)\;\rmd y_1\cdots\rmd y_{2q_0}\;,
$$
with, for all $(\xi_1,\xi_2)\in\mathbb{R}^{2}$,
\begin{equation}\label{e:exprg}
g(\xi_1,\xi_2)=
(1+|n\{\xi_1+\xi_2\}|)^{-2}
\frac{\left|\gamma_j\{\xi_1/\gamma_j\}\right|^{2(M-K)}\left|
\gamma_j\{\xi_{2}/\gamma_j \}\right|^{2(M-K)}}
{\left[(1+|\gamma_j\{\xi_1/\gamma_j \}|)
(1+|\gamma_j\{\xi_2/\gamma_j\}|)\right]^{2(M+\alpha)}}
\;.
\end{equation}
We now bound the integral $I_{n,j}$. Observe that for any
$y=(y_1,\cdots,y_{2q_0})\in (-n\gamma_j\pi,n\gamma_j\pi)^{2q_0}$
$$
|y_i+\dots+y_{q_0}|\leq n\gamma_j(q_0-i+1)\pi
\quad\text{and}\quad
|y_{q_0+i}+\dots+y_{2q_0}|\leq n\gamma_j(q_0-i+1)\pi\;.
$$
Thereafter, we set
\[
u_1=\frac{y_1+\cdots+y_{q_0}}{n},\,\cdots,u_{q_0}=\frac{y_{q_0}}{n},v_1=\frac{y_{q_0+1}+\cdots+y_{2q_0}}{n},\,\cdots,v_{q_0}=\frac{y_{2q_0}}{n}\;.
\]
Then
\begin{eqnarray*}
  I_{n,j}&\leq &
  c_0 \int_{u_1=-q_0\gamma_j\pi}^{q_0\gamma_j\pi}\int_{v_1=-q_0\gamma_j\pi}^{q_0\gamma_j\pi}
  g(u_1,v_1)\;J_{q_0,\gamma_j\pi}(u_1;2d\mathrm{1}_{q_0})
  J_{q_0,\gamma_j\pi}(v_1;2d\mathrm{1}_{q_0})
  \rmd u_1\rmd v_{1}\;,
\end{eqnarray*}
where we used the definition of $J_{m,a}(s;\beta)$ in Lemma~\ref{lem:Ja} with
the notation $\mathrm{1}_{q_0}$ for the $q_0$-dimensional vector with all
entries equal to 1, that is, we set $m=q_0$, $a=\gamma_j\pi$,
$\beta_1=\cdots=\beta_m=2d$ in~(\ref{e:Jma}).  We now apply
Lemma~\ref{lem:Ja}. Since $q_0<1/(1-2d)$, we are in Case (i) of
and we obtain
\begin{equation*}
J_{q_0,\gamma_j\pi}(z;2d\mathrm{1}_{q_0})
\leq C\;|z|^{-2\delta(m)}\;,\quad z\in\mathbb{R}\;,
\end{equation*}
for some constant $C>0$. This bound with the inequality $|\{u\}|\leq |u|$ and
the expression of $g$ given by~(\ref{e:exprg}) yields
\[
I_{n,j}\leq C\;\int_{-q_0\gamma_j\pi}^{q_0\gamma_j\pi}\int_{-q_0\gamma_j\pi}^{q_0\gamma_j\pi}\frac{
\left|\gamma_j\{u_1/\gamma_j\}\right|^{2(M-K-\delta(q_0))}
\left|\gamma_j\{v_{1}/\gamma_j\}\right|^{2(M-K-\delta(q_0))}\rmd u_1\rmd
v_{1}}{(1+n|\{u_1+v_{1}\}|)^2
\left[(1+\left|\gamma_j\{u_{1}/\gamma_j\}\right|)
(1+\left|\gamma_j\{v_{1}/\gamma_j\}\right|)\right]^{2(M+\alpha)}}\;.
\]
By $2\pi$--periodicity of $u\mapsto\{u\}$, we observe that the integrand is
$(2\pi\gamma_j)$-periodic with respect to both variables $u_1$ and $v_1$.
Thus the integral on $(-q_0\gamma_j\pi,q_0\gamma_j\pi)^2$
equals $q_0^2$ times the integral on $(-\gamma_j\pi,\gamma_j\pi)^2$. We get
that
\[
I_{n,j}\leq C\;
\int_{u_1=-\gamma_j\pi}^{\gamma_j\pi}\int_{v_1=-\gamma_j\pi}^{\gamma_j\pi}
\frac{\left|u_1\right|^{2(M-K-\delta(q_0))} \left|v_{1}\right|^{2(M-K-\delta(q_0))}\rmd
u_1\rmd
v_{1}}{(1+n|\{u_1+v_{1}\}|)^2\left(1+|u_1|\right)^{2(M+\alpha)}\left(1+|v_{1}|\right)^{2(M+\alpha)}}\;.
\]
By assumption $2(M-K-\delta(q_0))\geq0$, then for any $t\in\mathbb{R}$,
\[
\left|t\right|^{2(M-K-\delta(q_0))}\leq
\left(1+|t|\right)^{2(M-K-\delta(q_0))}\leq
\left(1+|t|\right)^{2(M-K)}\;.
\]
It implies that
\[
I_{n,j}\leq C\;
\int_{u_1=-\gamma_j\pi}^{\gamma_j\pi}\int_{v_1=-\gamma_j\pi}^{\gamma_j\pi}
\frac{\rmd u_1\rmd v_{1}}{(1+n|\{u_1+v_{1}\}|)^2 (1+|u_1|)^{2(K+\alpha)}(1+|v_{1}|)^{2(K+\alpha)}}\;.
\]
We now apply Lemma~\ref{lem:Jnbound} with
\[
S=2(K+\alpha),\,\beta_1=\beta_2=0\;.
\]
By assumption $S>1$. Then $I_{n,j}\leq C\;
n^{-1}$ and the conclusion follows from~(\ref{e:boundSnj1}).
\end{proof}

\section{The leading term of the scalogram and its asymptotic
  behavior}\label{s:leadingterms}
Suppose $q_0\geq 2$. We will show that the leading term of
$S_{n,j}$ is $S_{n,j}^{(q_0-1)}$ defined in~(\ref{e:Snjrrp1}). It
is an element of the chaos of order $2q_{0}-2(q_{0}-1)=2$ and
after renormalization it will converge to a Rosenblatt random
variable. We first study the asymptotic behavior of
$S_{n,j}-S_{n,j}^{(q_0-1)}$ which is a sum of random variables in
chaoses 4,6 up to $2q_{0}$. We actually show in the next result
that, under the normalization of $S_{n,j}^{(q_0-1)}$, this term is
negligible.

\begin{corollary}\label{cor:intSnj}
Assume $q_0\geq 2$ and $M\geq \delta(q_0)+K$. Then, as $j,n\to\infty$,
\begin{equation}\label{e:supTn0a}
n^{1-2d}\gamma_j^{-2(\delta(q_0)+K)} \left(\sum_{p=0}^{q_0-2}
p!{{q_{0}}\choose{p}}^2\|S_{n,j}^{(p)}\|_2\right)\rightarrow 0\;,
\end{equation}
\end{corollary}
\begin{proof}
The limit~(\ref{e:supTn0a}) is a direct consequence of
Propositions~\ref{pro:pneqzSD} and~\ref{pro:snjpnull}, observing that
$1-2d=1-2\delta(1)<1-2\delta(q_0-p)$ for all $p=1,2,\dots,q_0-2$ and that
$\delta(q_0)>0$ and $q_0\geq2$ imply $1-2d<1/2$.
\end{proof}

We consider the limit in distribution of the corresponding term
$n^{1-2d}\gamma_j^{-2(\delta(q_0)+K)}S_{n,j}^{(q_0-1)}$. With
Corollary~\ref{cor:intSnj}, this will provide the proof of Theorem~\ref{thm:LD}
in the case $q_0\geq2$. However, to cover the $m$-dimensional case with
$m\geq2$, we need to define a multivariate $S_{n,j}^{(p)}$ that will be denoted
by $\bS_{n,j}^{(p)}$. Let $0<p<q_0$.  Define a $\mathbb{C}^m$--valued function
$\widehat{\bkappa}_{j}^{(p)}$ by applying~(\ref{e:intrepKjp1}) component-wise
with $h_j$ replaced by $h_{\ell,j}$, $\ell=1,\dots,m$.  Define a
$\mathbb{C}^m$--valued function $\bg_p$ by~(\ref{e:bg}) with
$\widehat{\kappa}_{j}^{(p)}$ replaced by $\widehat{\bkappa}_{j}^{(p)}$.
Finally define $\bS_{n,j}^{(p)}$ as a $m$-dimensional random vector defined
by~(\ref{e:Snjrrp1}) with $g_p$ replaced by $\bg_p$.

\begin{proposition}\label{pro:lrd}
Suppose that Assumptions \textbf{A} hold with $2\leq q_0 <1/(1-2d)$ and
$M\geq K$. Then, for any diverging sequence $(n_j)$, as
$j\to\infty$, we have
\begin{equation}\label{e:norma}
n_j^{1-2d}\gamma_j^{-2(\delta(q_0)+K)} \bS_{n_j,j}^{(q_0-1)}
\overset{\mathcal{L}}{\longrightarrow}
f^*(0)^{q_0}\,\K_{q_0-1}
\, Z_d(1) \;.
\end{equation}
where $Z_{d}(1)$ and $\K_{q_0-1}$ are the same as in Theorem~\ref{thm:LD}.
% is the Rosenblatt process defined in~(\ref{e:harmros})is a finite positive
% constant defined in~(\ref{e:defKp}).
\end{proposition}
\begin{proof}
  Using~(\ref{e:Snjrrp}) component-wise with $p=q_0-1$, observing that
  $2q_0-2p=2$ and making the change of variable $y=n\gamma_j\xi$ in the
  multiple stochastic integral, we get, using the self-similarity of the Wiener
  process,
\begin{eqnarray}
\bS_{n,j}^{(q_0-1)} & =& \widehat{I}_{2}\left(
D_n\circ\Sigma_{2}(\gamma_{j}\times\cdot)\times
[\sqrt{f}\1_{(-\pi,\pi)}]^{\otimes 2}
\times \widehat{\bkappa}_{j}^{(q_0-1)}\right)\nonumber\\
& \stackrel{\tiny{d}}{=}&
\frac{1}{n\gamma_j}\;\widehat{I}_{2}\left(
D_n\circ\Sigma_{2}\left(n^{-1}\times\cdot\right)
\times\1_{(-\gamma_j\pi,\gamma_j\pi)}^{\otimes 2}
\left(n^{-1}\times\cdot\right)
\times \bof_{j}
\right)\;,\label{e:norma2}
\end{eqnarray}
where, for all $\xi\in\mathbb{R}^{2}$,
\begin{equation}\label{eq:fnj}
\bof_{j}(n\gamma_j \xi)=
\sqrt{f}^{\otimes 2}(\xi)
\times \widehat{\bkappa}_{j}^{(q_0-1)}(\xi)\;.
\end{equation}
Here $\stackrel{\tiny{d}}{=}$ means that the two vectors have same
distributions for all $n,j\geq 1$. We will use Lemma~\ref{lem:Dn} which
involves fractional parts. Let us express
$1_{(-\gamma_j\pi,\gamma_j\pi)}^{\otimes 2}$ as a sum of indicator
functions on the following pairwise disjoint domains,
\begin{equation*}
\Gamma_{j}^{(s)}=\{t=(t_1,t_2)\in
(-\gamma_{j}\pi,\gamma_{j}\pi)^{2},\,|t_1+t_2-2\pi s|<\pi\},
\quad s\in\mathbb{Z}\;.
\end{equation*}
Hence we obtain
\begin{equation}\label{eq:IjnDecomp}
\bS_{n,j}^{(q_0-1)} \stackrel{\tiny{d}}{=}
\frac{1}{n\gamma_j}
\sum_{s\in\mathbb{Z}}\bI_{n,j}^{(s)}\;.
\end{equation}
\begin{equation}\label{eq:Ijn}
\bI_{n,j}^{(s)}=
\widehat{I}_{2}\left(
D_n\circ\Sigma_{2}\left(n^{-1}\times\cdot\right)
\times\1_{\Gamma_{j}^{(s)}}\left(n^{-1}\times\cdot\right)
\times \bof_{j}\right)\;.
\end{equation}
Proposition~\ref{pro:lrd} follows from the following three
convergence results, valid for all fixed $m\in \mathbb{Z}$.
\begin{enumerate}[(a)]
\item If $s=0$, then, as $j\to\infty$,
\begin{equation}
  \label{eq:convL2s_fixed}
  (n_j\gamma_{j})^{-2 d}
\gamma_j^{-2(\delta(q_0-1)+K)}
\bI_{n_j,j}^{(0)}
  \stackrel{L^2(\Omega)}{\longrightarrow}
  (f^*(0))^{q_0} \; \K_{q_0-1} \;Z_{d}(1)  \;.
\end{equation}
\item We have, as $j\to\infty$,
  \begin{equation}
    \label{eq:convL2s_nz}
  \sup_{s\neq0} \mathbb{E}
\left[(n_j\gamma_{j})^{-4 d}\gamma_j^{-4(\delta(q_0-1)+K)}
\left|\bI_{n_j,j}^{(s)}\right|^2\right]\to 0\;.
  \end{equation}
\item We have, as $j\to\infty$,
  \begin{equation}
    \label{eq:convL2s_nzSum}
  \sum_{s\not\in \gamma_{j}\mathbb{Z}} \mathbb{E}
\left[(n_j\gamma_{j})^{-4 d}\gamma_j^{-4(\delta(q_0-1)+K)}
\left|\bI_{n_j,j}^{(s)}\right|^2\right]\to 0\;.
  \end{equation}
\end{enumerate}
To show that this is sufficient to prove the proposition, observe that, for any
$t=(t_1,t_2)\in\Gamma_{j}^{(s)}$, we have
$$
2\pi |s| - \pi < 2\pi |s| -
|t_1+t_2-2\pi s| \leq |t_1+t_2| < 2 \gamma_{j}\pi \;.
$$
Hence the domain $\Gamma_{j}^{(s)}$ is empty if $|s| > \gamma_j +1/2$. We
use~(\ref{eq:convL2s_nz}) for the two values $s=\gamma_j$ and $s=-\gamma_j$
and~(\ref{eq:convL2s_nzSum}) for the values
$s\notin\gamma_{j}\mathbb{Z}$. Thus~(\ref{eq:convL2s_nz})
and~(\ref{eq:convL2s_nzSum}) imply
$$
(n\gamma_{j})^{-2 d}\gamma_j^{-2(\delta(q_0-1)+K)}\sum_{s\neq0}\bI_{n,j}^{(s)}
  \stackrel{L^2(\Omega)}{\longrightarrow} 0 \;.
$$
Observe also that the normalizing factor in the left-hand side
of~(\ref{e:norma}) can be written as
\[
n^{1-2d}\gamma_j^{-2(\delta(q_0)+K)}=(n\gamma_j)\left((n\gamma_j)^{-2d}
%\gamma_j^{(q_0-1)(1-2d)-(2K+1)}
\gamma_j^{-2(\delta(q_0-1)+K)}
\right)\;,
\]
by using the definition of $\delta$ in~(\ref{e:ldparamq}). The last two
displays,~(\ref{eq:IjnDecomp}) and~(\ref{eq:convL2s_fixed})
yield~(\ref{e:norma}).

It only remains to prove~(\ref{eq:convL2s_fixed}),~(\ref{eq:convL2s_nz})
and~(\ref{eq:convL2s_nzSum}).

\indent a) We first show~(\ref{eq:convL2s_fixed}).  Since
$\bI_{n,j}^{(0)}$ and $Z_{d}(1)$ are defined as
stochastic integrals of order $2$, (\ref{eq:convL2s_fixed}) is
equivalent to the $L^2(\mathbb{R}^{2})$ convergence of the
normalized corresponding kernels. We show the
latter by a dominated convergence argument. These kernels are given
in~(\ref{eq:Ijn}) and~(\ref{e:harmros}) respectively. Observe that, as
$n\to\infty$, $D_n(\theta/n)\to(\rme^{\rmi\theta}-1)/(i\theta)$
by~(\ref{eq:dirichlet}), for all
$y\in\mathbb{R}^{2}$,
$$
D_n\left(n^{-1}(y_1+y_2)\right)\to
\frac{\exp(\rmi (y_1+y_2))-1}{(\rmi(y_1+y_2))}\;.
$$
By~(\ref{e:sdf0}), we have, as $(n\gamma_j)\to\infty$,
for all $y\in\mathbb{R}^{2}$,
$$
\sqrt{f}^{\otimes 2}(y/(n\gamma_j))\sim
f^*(0)\;(n\gamma_j)^{2d}\;|y_1|^{-d}|y_2|^{-d}\;.
$$
Now applying Lemma~\ref{lem:cvhjkp} to the $m$ entries of
$\widehat{\bkappa}_{j}^{(p)}$ with $p=q_0-1$, we get that,
as $j\to\infty$, for all $y\in\mathbb{R}^{2}$,
$$
\gamma_j^{(q_0-1)(1-2d)-(2K+1)}\widehat{\bkappa}_{j}^{(q_0-1)}(y/(n_j\gamma_j))\to
(f^*(0))^{q_0-1}\;\K_{q_0-1}\;.
$$
The last three convergences and $2\delta(q_0-1)=1-(q_0-1)(1-2d)$ yield the
pointwise convergence of the normalized kernels defining the
stochastic integrals appearing in the left-hand side
of~(\ref{eq:convL2s_fixed}) towards the kernel of the right-hand
side.

It remains to bound these kernels by an $L^2(\mathbb{R}^{2})$ function not
depending on $j,n$. We may take $m=1$ without loss of generality for this
purpose, since component-wise bounds are sufficient.  If
$y/n\in\Gamma_{j}^{(0)}$, we have, by Lemma~\ref{lem:Dn},
\begin{equation}
  \label{eq:DnBound_leadingterms}
\left|  D_n((y_1+y_2)/n)\right|\leq
C \; (1+|y_1+y_2|)^{-1} \;,
\end{equation}
for some constant $C>0$. By~(\ref{e:sdf0}), since $f^*$ is bounded,
we have, for all $y=(y_1,y_2)\in(-n\gamma_j\pi,n\gamma_j\pi)$
\begin{equation}
  \label{eq:fBound_leadingterms}
\left|(n\gamma_{j})^{-2 d}
\sqrt{f}^{\otimes 2}(y/(n\gamma_j))\right|\leq
C\,|y_1|^{-d}\;|y_2|^{-d}\;,
\end{equation}
where $C$ is a constant. Since $q_0-1<1/(1-2d)$, Lemma~\ref{lem:suphjkp}
implies that, for all $\zeta\in\mathbb{R}^2$ and some constant $C$,
\begin{equation}
  \label{eq:kappaBound_leadingterms}
  \left|\gamma_j^{-2(\delta(q_0-1)+K)} \widehat{\kappa}_{j}^{(q_0-1)}(\zeta)\right|
  \leq C\;.
\end{equation}
The bounds~(\ref{eq:DnBound_leadingterms}),~(\ref{eq:fBound_leadingterms})
and~(\ref{eq:kappaBound_leadingterms}) imply that
$(n\gamma_{j})^{-2
  d}\gamma_j^{-2(\delta(q_0-1)+K)}I_{n,j}^{(0)}=\widehat{I}_{2}(g)$ with
$$
|g(y)|^2 \leq C (1+|y_1+y_2|)^{-2} \;|y_1|^{-2d}|y_2|^{-2d},
\quad y=(y_1,y_2)\in\mathbb{R}^{2}\;,
$$
for some positive constant $C$.  Since we assumed $2<1/(1-2d)$. Then, applying
Lemma~\ref{lem:lem8.4CRTT} with $M_1=2$, $q=2$, and $a=0$, we obtain that this
function is integrable and the convergence~(\ref{eq:convL2s_fixed}) follows.

\indent b) Let us now prove~(\ref{eq:convL2s_nz}). Again we may take $m=1$
without loss of generality since the bound can be applied component-wise to
derive the case $m\geq2$.  Observe that the
bounds~(\ref{eq:fBound_leadingterms}) and~(\ref{eq:kappaBound_leadingterms})
can be used for $y/n\in\Gamma_{j}^{(s)}$, while the
bound~(\ref{eq:DnBound_leadingterms}) becomes
\begin{equation}
  \label{eq:DnBound_leadingterms2}
\left|  D_n(({y_1+y_2})/{n})\right|^2\leq
C \; (1+|y_1+y_2-2\pi n s |)^{-2} \;,
\end{equation}
Hence in this case, we obtain that
$(n\gamma_{j})^{-2
  d}\gamma_j^{-2(\delta(q_0-1)+K)}I_{n,j}^{(s)}=\widehat{I}_{2}(g)$ with
\begin{equation}\label{eq:fctg}
|g(y)|^2\leq
C (1+|y_1+y_2-2\pi n s|)^{-2} \;|y_1|^{-2d}\;|y_2|^{-2d},
\quad y=(y_1,y_2)\in\mathbb{R}^{2}\;,
\end{equation}
for some positive constant $C$. Using the assumption $2<1/(1-2d)$, from
Lemma~\ref{lem:lem8.4CRTT} applied with
$q=2$, $a=2\pi n s$ and $M_1=2$,  we get~(\ref{eq:convL2s_nz}).

\indent c) Finally we prove~(\ref{eq:convL2s_nzSum}) with $m=1$. We
need to further partition $\Gamma_{j}^s$ into
\begin{equation*}
\Gamma_{j}^{(s,\sigma)}=\{t\in\Gamma_{j}^s,\;
t/\gamma_j-2\pi \sigma\in (-\pi ,\pi)^2\},\quad
\sigma\in\mathbb{Z}^2\;.
\end{equation*}
Note that for all $t=(t_1,t_2)\in\Gamma_{j}^{(s,\sigma)}$, we have, for any
$i=1,2$,
$$
|2\pi\sigma_i|\leq |t_i/\gamma_j-2\pi\sigma_i| +  |t_i/\gamma_j|< 2\pi\;.
$$
Hence
$\Gamma_{j}^{(s,\sigma)}=\emptyset$ for all $\sigma$ out of the
integer rectangle $R=\{-1,0,1\}^2$. Then we obtain
$$
(n\gamma_{j})^{-2
  d}\gamma_j^{-2(\delta(q_0-1)+K)}I_{n,j}^{(s)}=
\sum_{\sigma\in R}\widehat{I}_{2}(g_\sigma^{(s)})\;,
$$
where, for all $y\in\mathbb{R}^{2}$,
$$
g_\sigma^{(s)}(y)=
(n\gamma_{j})^{-2 d}\gamma_j^{-2(\delta(q_0-1)+K)}
D_n\circ\Sigma_2(y/n)
\times\1_{\Gamma_{j}^{(s,\sigma)}}(y/n)
\times f_{j}(y)\;.
$$
Since $R$ is a finite set, to obtain the limit~(\ref{eq:convL2s_nzSum}), it is
sufficient to show that, for any fixed $\sigma\in R$, as $j,n\to\infty$,
  \begin{equation}
    \label{eq:convL2s_nzSumSigma_fixed}
  \sum_{s\not\in \gamma_{j}\mathbb{Z}}
\int \left|g_\sigma^{(s)}(y)\right|^2\;\rmd^{2} y
\to 0\;.
  \end{equation}
  For $\zeta\in 2\pi \sigma+(-\pi ,\pi)^2$, we use a sharper bound
  than~(\ref{eq:kappaBound_leadingterms}), namely, by Lemma~\ref{lem:suphjkp},
\begin{equation}\label{e:boundkap61}
\left|\gamma_j^{-2(\delta(q_0-1)+K)}\widehat{\kappa}_j^{(q_0-1)}(\zeta)\right|^2\leq
 C\;k_j^{\otimes 2}(\zeta-2\pi\sigma)\quad\text{where}\quad k_j(u)=
(1+\gamma_j|u|)^{-2\delta(q_0-1)}\;.
\end{equation}
With~(\ref{eq:fBound_leadingterms})
and~(\ref{eq:DnBound_leadingterms2}), it follows that
\begin{equation}
  \label{eq:g_sig_s_bound}
  \left|g_\sigma^{(s)}(y)\right|^2\leq
  C \; \frac{k_j^{\otimes 2}(y/(n\gamma_j)-2\pi\sigma)}
  {(1+|y_1+y_2-2\pi n s|)^{2}} \;|y_1|^{-2d}\;|y_2|^{-2d},
  \quad y=(y_1,y_2)\in\mathbb{R}^{2}\;.
\end{equation}
Let us set $w=(w_1,w_2)$ with $w_1=y_1/(n\gamma_j)-2\pi\sigma_1$ and
$w_2=y_2/(n\gamma_j)-2\pi\sigma_2$.
Using the bound~(\ref{eq:g_sig_s_bound}) and that
$y/n\in\Gamma_{j}^{(s,\sigma)}$ implies $w\in
\Delta_{j}^{(s,\sigma)}$ with
$$
\Delta_{j}^{(s,\sigma)}=\{(w_1,w_2)\in
(-\pi,\pi)^2,\,|\gamma_j(w_1+w_2)-2\pi(s-\gamma_j(\sigma_1+\sigma_2))|<\pi\}\;,
$$
we get
$$
\int \left|g_\sigma^{(s)}(y)\right|^2\;\rmd^{2} y \leq C
(n\gamma_j)^{2(1-2d)}\int_{\Delta_{j}^{(s,\sigma)}}
\frac{k_j^{\otimes 2}(w)\;|w_1+2\pi\sigma_1|^{-2d}|w_2+2\pi\sigma_2|^{-2d}}
{(1+n|\gamma_j(w_1+w_2)-2\pi(s-\gamma_j(\sigma_1+\sigma_2))|)^{2}}\;\rmd^2
w \;,
$$
Since $|w_i\pm 2\pi|>\pi>|w_i|$ for $w\in\Delta_{j}^{(s,\sigma)}$, we
have for $\sigma\in R$,
\begin{equation}\label{eq:int_g_sigma_s_Delta}
\int
\left|g_\sigma^{(s)}(y)\right|^2\;\rmd^{2} y \leq C
(n\gamma_j)^{2(1-2d)}\int_{\Delta_{j}^{(s,\sigma)}}
\frac{k_j^{\otimes 2}(w)\;|w_1|^{-2d}|w_2|^{-2d}}
{(1+n|\gamma_j(w_1+w_2)-2\pi
(s-\gamma_j(\sigma_1+\sigma_2))|)^{2}} \;\rmd^2 w \;.
\end{equation}
We shall apply Lemma~\ref{lem:lem8.4CRTT} after having conveniently bounded
$k_j$ in the numerator of the previous ratio. Let $\beta<1$ to be set later
arbitrarily close to 1.  Since $2\delta(q_0-1)\geq \beta-2d+2\delta(q_0)$, we
have
\begin{align*}
k_j(u)&=(1+\gamma_j|u|)^{-2\delta(q_0-1)}\\
&\leq
(1+\gamma_j|u|)^{2d-\beta}(1+\gamma_j|u|)^{-2\delta(q_0)}\;.
\end{align*}
Observe that, for all $w\in\Delta_{j}^{(s,\sigma)}$ we have
$$
\gamma_j (|w_1|\vee |w_2|)\geq
\gamma_j |(w_1+w_2)/2|\geq
\pi(|s-\gamma_j(\sigma_1+\sigma_2)|-1/2) \geq
\pi|s-\gamma_j(\sigma_1+\sigma_2)|/2\;.
$$
In the last inequality, we used that $s\not\in\gamma_j\mathbb{Z}$ and that $s$,
$\gamma_j$, $\sigma_1$ and $\sigma_2$ are integers so that
$|s-\gamma_j(\sigma_1+\sigma_2)|\geq 1$.

Using $0< q_0<1/(1-2d)$, we have $2\delta(q_0)>0$, and, choosing $\beta$ close
enough to 1, we have $\beta-2d>0$. Hence, the last two displays yield, for all
$w\in\Delta_{j}^{(s,\sigma)}$ with $s\not\in\gamma_j\mathbb{Z}$,
\begin{equation}\label{e:k61b}
k_j^{\otimes2}(w)
\leq|\gamma_jw_1|^{2d-\beta}|\gamma_jw_2|^{2d-\beta}
(1+\pi|s-\gamma_j(\sigma_1+\sigma_2)|/2)^{-2\delta(q_0)} \;.
\end{equation}
Inserting this bound in~(\ref{eq:int_g_sigma_s_Delta}) and setting $t=n\gamma_j
w$,
 we
obtain
\begin{eqnarray*}
&&\int \left|g_\sigma^{(s)}(y)\right|^2\;\rmd^{2} y\\
&&\leq C \frac{n^{-4d+2\beta}}
{|s-\gamma_j(\sigma_1+\sigma_2)|^{2\delta(q_0)} }
\int_{\mathbb{R}^2}
\frac{|t_1t_2|^{-\beta}}{(1+|t_1+t_2-2\pi n(s-\gamma_j(\sigma_1+\sigma_2))|)^{2}}\;\rmd^2 t\;.
\end{eqnarray*}
For $\beta$ close enough to 1, we may apply Lemma~\ref{lem:lem8.4CRTT} with
$q=2$, $d=\beta/2$, $M_1=2$ and $a=2\pi n(s-\gamma_j(\sigma_1+\sigma_2))$ to
bound the previous integral.
Using again that $s\not\in\gamma_j\mathbb{Z}$ and that $s$, $\gamma_j$,
$\sigma_1$ and $\sigma_2$ are integers, we have
$|a|\geq 2 \pi n$ and thus $1+|a|\asymp |a|$.
We get, for all $s\not\in\gamma_j\mathbb{Z}$
$$
\int \left|g_\sigma^{(s)}(y)\right|^2\;\rmd^{2} y
\leq C n^{1-4d}\;|s-\gamma_j(\sigma_1+\sigma_2)|^{1-2\delta(q_0)-2\beta}  \;,
$$
where $C$ is some positive constant.

Now choose $\beta$ close enough to 1
so that $2\delta(q_0)+2\beta-1>1$. It follows that
$$
\sum_{k\neq0}|k|^{1-2\delta(q_0)-2\beta}<\infty\;.
$$
Since our assumptions imply $d>1/4$, the last two displays
imply~(\ref{eq:convL2s_nzSumSigma_fixed}) and the proof is finished.
\end{proof}
\section{Proof of the main results}\label{s:proofs-main}
\begin{proof}[Proof Theorem~\ref{thm:LD}]
We first prove the result in Case~\ref{item:thm-case-gauss}.  In this case
$q_0=1$ and thus $H_{q_0}(X_t)=X_t$.  Let $(v(s))_{s\in\mathbb{Z}}$ be the
Fourier coefficients of $\sqrt{2\pi f}$, so that the convergence
$$
\sqrt{2\pi\,f(\lambda)}=\widehat{v}(\lambda)
=\sum_{s\in\mathbb{Z}}v(s)\rme^{-\rmi\lambda s}
$$
holds in $L^2(-\pi,\pi)$. It follows that $\{X_t\}_{t\in\mathbb{Z}}$ can be
represented as
$$
X_t=\sum_{s\in\mathbb{Z}}v(t-s)\xi_s,\quad t\in\mathbb{Z}\;,
$$
where $\{\xi_t\}_{t\in\mathbb{Z}}$ is an i.i.d. sequence of standard Gaussian
r.v.'s. Applying~(\ref{e:W2}) with $H_{q_0}(X_t)=X_t$ we obtain that
\begin{equation}
  \label{eq:WZ}
\bW_{j,k}= \gamma_j^{d+K}
\left[
\begin{array}{c}
Z_{1,j,k}\\
\vdots\\
Z_{m,j,k}
\end{array}
\right]  \;,
\end{equation}
where
$$
Z_{\ell,j,k} = \sum_{t\in\mathbb{Z}}v_{\ell,j}(\gamma_j k -t)\xi_t \,
$$
with
$$
v_{\ell,j}(u) = \gamma_j^{-d-K}
\sum_{s\in\mathbb{Z}} h_{\ell,j}^{(K)}(u-s)\;v(s),\quad
u\in\mathbb{Z}\;.
$$
Hence
$$
\widehat{v}_{\ell,j}(\lambda) = \gamma_j^{-d-K}\widehat{h}_{\ell,j}^{(K)}(\lambda)
\widehat{v}(\lambda) =\gamma_j^{-d-K}\;\sqrt{2\pi\,f(\lambda)} \;
\widehat{h}_{\ell,j}^{(K)}(\lambda),\quad\lambda\in(-\pi,\pi)\;.
$$
Observe that~(\ref{e:sdf0}),~(\ref{e:majoHj}) and~(\ref{e:HjkVSHj}) imply,
for some positive constant $C$,
$$
\left|\widehat{v}_{\ell,j}(\lambda)\right|\leq
C\gamma_j^{1/2}
\frac{|\gamma_j\lambda|^{M-(K+d)}}{(1+\gamma_j|\lambda|)^{\alpha+M}},
\quad\lambda\in(-\pi,\pi)\;.
$$
On the other hand, ~(\ref{e:sdf0}),~(\ref{EqLimHj}) and~(\ref{e:HjkVSHj}) imply
$$
\lim_{j \to +\infty}\gamma_j^{-1/2}\widehat{v}_{\ell,j}(\gamma_j^{-1}\lambda)
\rme^{\rmi\Phi_j(\lambda)}
=\sqrt{2\pi f^*(0)}|\lambda|^{-(K+d)}
\widehat{h}_{\ell,\infty}(\lambda),\quad\lambda\in\mathbb{R},\;\ell=1,\dots,m\;.
$$
Thus, if $M\geq K+d$, Assumptions \textbf{A} imply Condition \textbf{B} in
\cite{roueff-taqqu-2009} with $N=m$, $\delta=\alpha+K+d$, $\lambda_{i,j}=
\lambda_{i,\infty}=0$, $\Phi_{i,j}=\Phi_j$,
$v_{i,j}^*=(2\pi)^{-1/2}\widehat{v}_{i,j}$ and
$v_{i,\infty}^*(\lambda)=\sqrt{f^*(0)}|\lambda|^{-(K+d)}
\widehat{h}_{i,\infty}(\lambda)$ for $i=1,\dots,N$ and $j\geq1$.
Moreover we may apply Theorem~1 in \cite{roueff-taqqu-2009}
and obtain, as $j\to\infty$,
$$
n_j^{-1/2}\sum_{k=0}^{n_j-1}
\left[
\begin{array}{c}
Z_{1,j,k}^2 -\mathbb{E}[Z_{1,j,k}^2]\\
\vdots\\
Z_{N,j,k}^2 -\mathbb{E}[Z_{N,j,k}^2]
\end{array}
\right]
\overset{\mathcal{L}}{\longrightarrow}
\mathcal{N}(0,\Gamma) \; ,
$$
where $\Gamma$ is the $m\times m$ covariance matrix defined
by~(\ref{eq:GammaDef}). Since, by~(\ref{e:snjm}) and~(\ref{eq:WZ}),
$n_j^{1/2}\gamma_j^{-2(d+K)} \overline{\bS}_{n_j,j}$ is the left-hand side of the
last display, we get~(\ref{eq:CenteredZ}).

We now consider Case~\ref{item:thm-case-rozen}.
Applying
the basic decomposition~(\ref{e:Snjdec}) to each entries of
$\overline{\bS}_{n,j}$, Corollary~\ref{cor:intSnj}  and
Proposition~\ref{pro:lrd} show that the leading term is obtained for
$p=q_0-1$. Moreover the latter proposition specifies the limit.
\end{proof}

\begin{proof}[Proof of Theorem~\ref{thm:multiscale}]
We first
prove~(\ref{eq:wav-sp-asymp}). Applying~(\ref{eq:wav-sp_def}),~(\ref{eq:Wjk_representation})
(with $h_j$ replaced by $g_j$) and the isometry
property~(\ref{e:cov-multiple}), we have
$$
{\sigma}^2_{j} = q_0! \int_{(-\pi,\pi)^{q_0}}
\frac{\left|\widehat{g}_j\circ\Sigma_{q_0}(\xi)\right|^2}
{|1-\rme^{-\rmi \Sigma_{q_0}(\xi)}|^{2K}} f^{\otimes q_0}(\xi)\,\rmd^{q_0}\xi\;.
$$
Setting $\xi=2^{-j}\lambda$, we get
$$
{\sigma}^2_{j} = q_0!2^{-j(q_0-1)} \int_{(-2^j\pi,2^j\pi)^{q_0}}
\frac{\left|2^{-j/2}\widehat{g}_j\circ\Sigma_{q_0}(2^{-j}\lambda)\right|^2}
{|1-\rme^{-\rmi \Sigma_{q_0}(2^{-j}\lambda)}|^{2K}} f^{\otimes
  q_0}(2^{-j}\lambda)\,\rmd^{q_0}\lambda\;.
$$
Using Assumption (W-\ref{ass:w-b}) on $g_j$,
and Condition~(\ref{e:sdf}) with $f^\ast$ bounded, the integrand is
bounded, up to a multiplicative constant, by
$$
2^{2j(K+d q_0)}
(1+|\Sigma_{q_0}(\lambda)|)^{-2(\alpha+K)}
\prod_{i=1}^{q_0}|\lambda_i|^{-2d}\,,
$$
since $(|x|/(1+|x|))^M\leq (|x|/(1+|x|))^K$ and $|1-\rme^{-\rmi x}|\asymp|x|$.
The displayed bound is integrable by Lemma~\ref{lem:lem8.4CRTT} with
$M_1=2(\alpha+K)$, $q=q_0$ and $a=0$. By dominated convergence, Assumption
(W-\ref{ass:w-c}) on $(g_j)$ and continuity of $f^*$ at zero, we get, as
$j\to\infty$,
$$
2^{-2j(K+d q_0-(q_0-1)/2)}{\sigma}^2_{j} \to
q_0! \, (f^*(0))^{q_0}\,
 \int_{\mathbb{R}^{q_0}}
\frac{|\widehat{g}_\infty\circ\Sigma_{q_0}(\lambda)|^2}
{|\Sigma_{q_0}(\lambda)|^{2K}}\,
\prod_{i=1}^{q_0}|\lambda_i|^{-2d}\,\rmd^{q_0}\lambda\;.
$$
Using~(\ref{e:ldparamq}) and the definition $L_{q_0}(\widehat{g}_\infty)$ in~(\ref{e:defKp}), we
obtain~(\ref{eq:wav-sp-asymp}).

To prove the convergence of the scalogram, we shall apply
Theorem~\ref{thm:LD}(b) with a sequence of multivariate filters
$(\bh_j)_{j\geq0}$.  To illustrate how this is done, suppose, for example, that
we want to study the joint behavior of $W_{j-u,k}$ for $u\in\{0,1\}$. Recall
that $j-1$ is a finer scale than $j$. Following the framework of
\cite{roueff-taqqu-2009b}, we consider the multivariate coefficients
$\bW_{j,k}=(W_{j,k},\,W_{j-1,2k},\,W_{j-1,2k+1})$, since, in addition to the
wavelet coefficients $W_{j,k}$ at scale $j$, there are twice as many wavelet
coefficients $W_{j-1,2k},\,W_{j-1,2k+1}$ at scale $j-1$.  These coefficients
can be viewed as the output of a multidimensional filter $\bh_j$ defined as
$\bh_j(\tau)=(h_j(\tau),\,h_{j-1}(\tau),\,h_{j-1}(\tau+2^{j-1}))$.
These three entries correspond to $(u,v)$ equal to $(0,0)$, $(1,0)$ and
$(1,1)$, respectively, in the general case below.

In the general case, each $\bh_j$ is defined as follows.
For all, $j\geq0$, $u\in\{0,\dots,j\}$
and $v\in\{0,\dots,2^u-1\}$, let $\ell=2^u+v$ and define a filter $h_{\ell,j}$
by
\begin{equation}
  \label{eq:multiscale-filers}
  h_{\ell,j}(t)=g_{j-u}(t+2^{j-u} v),\quad t\in\mathbb{Z}\;.
\end{equation}
Applying this definition and~(\ref{eq:wav_coeff_def1_dyad}) with $\gamma_j=2^j$, we
get
$$
W_{j-u,2^uk+v}=\sum_{t\in\mathbb{Z}}h_{\ell,j}(2^j k-t)Y_t\;.
$$
These coefficients are stored in a vector $\bW_{j,k}=[W_{\ell,j,k}]_\ell$, say
of length $m=2^p-1$,
\begin{equation}
  \label{eq:Wljk-multiscale}
W_{\ell,j,k}=W_{j-u,2^uk+v}, \quad\ell=2^u+v=1,2,\dots,m\,,
\end{equation}
which corresponds to the multivariate wavelet coefficient~(\ref{e:bW1}) with
$\bh_{j}(t)$ having components $h_{\ell,j}(t)$, $\ell=1,2,\dots,m$ defined
by~(\ref{eq:multiscale-filers}). This way of proceeding allows us to express
the vector $[\widehat{\sigma}^2_{j-u}-{\sigma}^2_{j-u}]_{u=0,\dots,p-1}$ as a
linear function, up to a negligible remainder, of the vector
$\overline{\bS}_{n_j,j}$
defined by~(\ref{e:snjm}).  Indeed observe that~(\ref{eq:nj_def}) implies, for
any fixed $u$
\begin{equation}
  \label{eq:njVSnj-u}
n_{j-u}=2^u n_j+O(1) \;.
\end{equation}
Hence~(\ref{eq:scalo_def}) and~(\ref{eq:wav-sp_def})
imply, for any fixed $u$,
\begin{align*}
\widehat{\sigma}^2_{j-u}-{\sigma}^2_{j-u}&=
\frac1{n_{j-u}}\sum_{k=0}^{n_{j-u}}
\left(W_{j-u,k}^2-\mathbb{E}[W_{j-u,k}^2]\right)\\
&=\frac1{n_{j-u}}\sum_{k=0}^{2^u n_{j}-1}
\left(W_{j-u,k}^2-\mathbb{E}[W_{j-u,k}^2]\right)
+O_P(\sigma_{j-u}^2/n_{j-u}),\quad
j\geq u\;.
\end{align*}
Expanding $\sum_{k=0}^{2^u n_{j}-1}$ as
$\sum_{v=0}^{2^u-1}\sum_{k'=0}^{n_{j}-1}$ with $k=k'2^{u}+v$ and
applying~(\ref{eq:Wljk-multiscale}) and the last display, we obtain, for all
$j\geq p$,
\begin{align}
\nonumber
\widehat{\sigma}^2_{j-u}-{\sigma}^2_{j-u}
&=\frac{1}{n_{j-u}}\sum_{v=0}^{2^u-1}\sum_{k'=0}^{n_{j}-1}
\left(W_{2^u+v,j,k}^2-\mathbb{E}[W_{2^u+v,j,k}^2]\right)
+O_P(\sigma_{j-u}^2/n_{j-u})\\
  \label{eq:multisvar-to-multiscale}
&=\frac{n_j}{n_{j-u}}\sum_{v=0}^{2^u-1}
\overline{S}_{n_j,j}(2^u+v)+O_P(\sigma_{j-u}^2/n_{j-u}),
\quad u=0,\dots,p-1\;,
\end{align}
where we denoted the entries of
$\overline{\bS}_{n_j,j}$ in~(\ref{e:snjm}) as
$[\overline{S}_{n_j,j}(\ell)]_{\ell=1,\dots,m}$.

Let us now check that $(\bh_j)$ satisfies the assumptions of
Theorem~\ref{thm:LD}. By hypothesis $\{g_j\}$ verifies Assumptions
(W-\ref{ass:w-a})--(W-\ref{ass:w-c}). Hence, by~(\ref{eq:multiscale-filers}),
$\{\bh_j\}$ satisfies (W-\ref{ass:w-a}). We further have that, for $\ell=2^u+v$
with $u\in\{0,\dots,p-1\}$ and $v\in\{0,\dots,2^u-1\}$,
$$
\widehat{h}_{\ell,j}(\lambda)=\widehat{g}_{j-u}(\lambda)
\rme^{\rmi 2^{j-u}v\lambda},\,\lambda\in(-\pi,\pi)\;.
$$
Hence (W-\ref{ass:w-b}) follows from the assumption on $g_j$.  Using that
$\gamma_j=2^j$, Condition (W-\ref{ass:w-c}) also follows with $\Phi_j\equiv0$
and
\begin{equation}
  \label{eq:multiscale_infty_filter}
\widehat{h}_{\ell,\infty}(\lambda)=
2^{-u/2}\widehat{g}_{\infty}(2^{-u}\lambda)\rme^{\rmi2^{-u}v\lambda}\;.
\end{equation}
We can thus apply Theorem~\ref{thm:LD} and obtain~(\ref{eq:caseBasymp_dist}),
that is,
$$
n_j^{1-2d} 2^{-2j(\delta(q_{0})+K)}
\overline{\bS}_{n_j,j}
\overset{\mathcal{L}}{\longrightarrow}
f^*(0)^{q_0}\,\K_{q_0-1}\,  Z_d(1) \;.
$$
with $\K_{q_0-1}=[L_{q_0-1}(\widehat{h}_{\ell,\infty})]_{\ell=1,\dots,m}$.
By~(\ref{eq:multiscale_infty_filter}) and~(\ref{e:defKp}), it turns out that,
for $\ell=2^u+v$ with $u\in\{0,\dots,p-1\}$ and $v\in\{0,\dots,2^u-1\}$,
\begin{align*}
L_{q_0-1}(\widehat{h}_{\ell,\infty})
&=\int_{\mathbb{R}^{q_0-1}}\frac{|2^{-u/2}\widehat{g}_\infty(2^{-u}(t_1+\cdots+t_{p}))|^{2}}
{|u_1+\cdots+u_{q_0-1}|^{2K}}\; \prod_{i=1}^{q_0-1} |t_i|^{-2d}\; \rmd
t_1\cdots\rmd t_{q_0-1} \\
&=2^{-u-2Ku-2d(q_0-1)u+u(q_0-1)}
L_{q_0-1}(\widehat{g}_{\infty}) \\
&=2^{-2u(\delta(q_{0})+K)+u(2d-1)}
L_{q_0-1}(\widehat{g}_{\infty}) \;,
\end{align*}
after the change of variables $s_i=2^{-u}t_i$, $i=1,\dots,q_0-1$  and the
definition of $\delta(q_0)$ in~(\ref{e:ldparamq}).
Using the last two displays, we obtain
that, as $j\to\infty$,
$$
\left\{n_j^{1-2d} 2^{-2(j-u)(\delta(q_{0})+K)}
\overline{\bS}_{n_j,j}(2^u+v)\right\}_{u,v}
\overset{\mathrm{fidi}}{\longrightarrow}
\left\{2^{u(2d-1)}\,L_{q_0-1}(\widehat{g}_{\infty})
f^*(0)^{q_0}\,  Z_d(1)\right\}_{u,v}\;,
$$
where $(u,v)$ take values $u=0,\dots,p-1$ and $v=0,\dots,2^u-1$.  Note that the
right-hand side does not depend on $v$.  By~(\ref{eq:nj_def}), we have
$n_j/n_{j-u}\sim2^{-u}$ and by~(\ref{eq:multisvar-to-multiscale}), we have
${\sigma}^2_{j-u} \sim q_0!  \, (f^*(0))^{q_0} \,
L_{q_0}(\widehat{g}_{\infty})\; 2^{2 (j-u) (\delta(q_0)+K)}$. Thus the last
display yields
$$
\left\{n_j^{1-2d}
\frac1{{\sigma}^2_{j-u}}
\frac{n_j}{n_{j-u}}\sum_{v=0}^{2^u-1}\overline{\bS}_{n_j,j}(2^u+v)\right\}_u
% &=
% \frac{n_j}{n_{j-u}}
% 2^{2u(\delta(q_{0})+K)}
% \sum_{v=0}^{2^u-1}
\overset{\mathrm{fidi}}{\longrightarrow}
\left\{
2^{u(2d-1)}\,
\frac{L_{q_0-1}(\widehat{g}_{\infty})}{q_0! L_{q_0}(\widehat{g}_{\infty})}\,
  Z_d(1)\right\}_u \;,
$$
where $u=0,\dots,p-1$.
Applying~(\ref{eq:multisvar-to-multiscale}), we have
\begin{align*}
n_j^{1-2d}
     \left(\frac{\widehat{\sigma}^2_{j-u}}{{\sigma}^2_{j-u}}-1\right)&=
n_j^{1-2d}\,\frac1{{\sigma}^2_{j-u}}\,
     \left(\widehat{\sigma}^2_{j-u}-{\sigma}^2_{j-u}\right)\\
  &=n_j^{1-2d}\,\frac1{{\sigma}^2_{j-u}}\,
\frac{n_j}{n_{j-u}}\sum_{v=0}^{2^u-1}
\overline{S}_{n_j,j}(2^u+v)+O_P(n_j^{1-2d}/n_{j-u})\;.
\end{align*}
By~(\ref{eq:njVSnj-u}), $n_j^{1-2d}/n_{j-u}\sim2^u n_j^{-2d}\to0$ since $u$ is
constant. Hence~(\ref{eq:scalo-multiscale-asymp}) follows from the last two
displays.
\end{proof}

\appendix

\section{Technical lemmas}\label{sec:techlemma}
\subsection{Asymptotic behavior of the kernel $\widehat{\kappa}_{j}^{(p)}$}\label{sec:asymptkappa}

The following result provides a bound of
$\widehat{\kappa}_{j}^{(p)}$ defined in~(\ref{e:intrepKjp1}), in
the case where $p>0$. It is used in the proof of
Proposition~\ref{pro:pneqzSD}.

\begin{lemma}\label{lem:suphjkp}
  Suppose that Assumptions \textbf{A} hold with $m=1$ and $M\geq K$, and let
  $0<p<1/(1-2d)$. Then there exists some $C_1>0$ such that for all
  $(\xi_1,\xi_2)\in \mathbb{R}^2$ and $j\geq0$,
\begin{equation}\label{e:majohjkp}
|\widehat{\kappa}_{j}^{(p)}(\xi_1,\xi_2)| \leq C_1
\frac{\gamma_j^{2(\delta(p)+K)}}{(1+\gamma_j|\{\xi_1\}|)^{\delta(p)}
(1+\gamma_j|\{\xi_2\}|)^{\delta(p)}}\;.
\end{equation}
\end{lemma}
\begin{proof}
  By $(2\pi)$-periodicity of $\widehat{\kappa}_{j}^{(p)}(\xi_1,\xi_2)$ along
  both variables $\xi_1$ and $\xi_2$, we may take
  $\xi_1,\xi_2\in[-\pi,\pi]$. Set for all $i\in\{1,\cdots,p\}$,
\[
\mu_i=\gamma_j\left(\lambda_i+\cdots+\lambda_p\right),
\]
in the integral~(\ref{e:intrepKjp}). Then by~(\ref{e:sdf})
and~(\ref{eq:M_replaced_by_K}), there exists a constant $C$ independent of
$j$ such that for all $(\xi_1,\xi_2)\in[-\pi,\pi]^2$,
\begin{eqnarray*}
|\widehat{\kappa}_j^{(p)}(\xi_1,-\xi_2)| \leq
C\|f^*\|^p_{\infty}\gamma_j^{2K+2\delta(p)}\int_{-\gamma_j
p\pi}^{\gamma_j p\pi}\frac{J_{p,\gamma_j \pi}(\mu_1;2d)\rmd \mu_1}
{\prod_{i=1}^2
\left(1+\gamma_j\left|\left\{{\mu_1}/{\gamma_j}+\xi_i\right\}\right|\right)^{K+\alpha}}\;,
\end{eqnarray*}
where $J_{p,a}$ is defined in Lemma~\ref{lem:Ja}. Applying Lemma~\ref{lem:Ja} ($\beta=2d, a=\gamma_j\pi$),
there exists some constant $C>0$ depending only on $p,d$ such that
for any $\mu_1\in \mathbb{R}^*$,
\begin{align}\label{e:boundJp}
J_{p,\gamma_j \pi}(\mu_1,2d)&\leq C
|\mu_1|^{-(p(1-2d)-1)}=C|\mu_1|^{-2\delta(p)}\;.
\end{align}
Hence there exists $C_1>0$ such that, for all $(\xi_1,\xi_2)\in
[-\pi,\pi]^2$,
\[
|\widehat{\kappa}_j^{(p)}(\xi_1,-\xi_2)| \leq C_1
\gamma_j^{2K+2\delta(p)}\int_{-p\gamma_j
\pi}^{p\gamma_j \pi}\frac{|\mu_1|^{-2\delta(p)}\rmd \mu_1}
{\prod_{i=1}^2\left(1+\gamma_j\left|\left\{{\mu_1}/{\gamma_j}+\xi_i\right\}\right|\right)^{K+\alpha}}\;.
\]
Using the Cauchy--Schwartz
inequality yields
\begin{equation}\label{e:boundK1}
|\widehat{\kappa}_j^{(p)}(\xi_1,-\xi_2)| \leq C_1 \gamma_j^{2(K+\delta(p))}\prod_{i=1}^{2}\left(\int_{-p\gamma_j
\pi}^{p\gamma_j \pi}\frac{|\mu_1|^{-2\delta(p)}\rmd \mu_1}
{(1+\left|\gamma_j\left\{{\mu_1}/{\gamma_j}+\xi_i\right\}\right|)^{2(K+\alpha)}}\right)^{1/2}\;.
\end{equation}

We now use that
$$
\int_{-p\gamma_j
\pi}^{p\gamma_j \pi}
\frac{|\mu_1|^{-2\delta(p)}\,\rmd
\mu_1 }
{(1+\left|\gamma_j\left\{\mu_1/{\gamma_j}+\xi\right\}\right|)^{2(K+\alpha)}}
\leq \sum_{|s|<(p+1)/2} \int_{I(s)} \frac{|\mu_1|^{-2\delta(p)}\,\rmd
\mu_1 }
{(1+\left|\mu_1+\gamma_j(\xi-2\pi s)\right|)^{2(K+\alpha)}}\;,
$$
where $I(s)$ denotes the interval $-\gamma_j\xi+2\pi
s\gamma_j+[-\gamma_j\pi,\gamma_j\pi]$.  Since we have here supposed that
$\delta(p)>0$, we may apply Lemma~\ref{lem:lem8.4CRTT} with $d=\delta(p)$,
$q=1$, $a=-\gamma_j(\xi-2\pi s)$ and $M_1=2(K+\alpha)$. We get
$$
\int_{-p\gamma_j
\pi}^{p\gamma_j \pi}
\frac{|\mu_1|^{-2\delta(p)}\,\rmd
\mu_1 }
{(1+\left|\gamma_j\left\{\mu_1/{\gamma_j}+\xi\right\}\right|)^{2(K+\alpha)}}
\leq C\, \sum_{|s|<(p+1)/2}(1+\gamma_j|\xi-2\pi s|)^{-2\delta(p)}\;,
$$
for some positive constant $C$. Since $|\xi|\leq\pi$, we have, for any non-zero
integer $s$, $|\xi-2\pi s|\geq (2|s|-1)\pi\geq \pi\geq |\xi|$. Hence all the
terms in the last sum are at most equal to the term corresponding to $s=0$.
This, with~(\ref{e:boundK1}), yields~(\ref{e:majohjkp}).
\end{proof}

Next we derive the limit of $\widehat{\kappa}_{j}^{(p)}$, rescaled and normalized, as $j\to\infty$. The result is
used in the proof of Proposition~\ref{pro:lrd}.
\begin{lemma}\label{lem:cvhjkp}
  Suppose that Assumptions \textbf{A} hold with $m=1$ and $M\geq K$, and let
  $0<p<1/(1-2d)$. Let $(z_{j})_{j\geq1}$ be a sequence in
  $\mathbb{R}^{2}$ converging to the origin. Then, as $j\to\infty$,
\[
\gamma_j^{p(1-2d)-(2K+1)}\widehat{\kappa}_{j}^{(p)}\left(z_{j}/\gamma_j\right)
\to (f^*(0))^{p}\;L_p(\widehat{h}_{\infty})\;,
\]
where $L_p(\widehat{h}_{\infty})$ is the finite positive constant defined
by~(\ref{e:defKp}).
\end{lemma}
\begin{proof}
  From~(\ref{EqMajoHinf})and~(\ref{e:defKp}) with $M\geq K$ we get that
  $|\widehat{h}_{\infty}(\lambda)|/|\lambda|^K\leq (1+|\lambda|)^{-\alpha-K}$.
  The fact that $L_p(\widehat{h}_{\infty})<\infty$ follows from
  Lemma~\ref{lem:lem8.4CRTT} applied with $a=0$, $p=q$ and $M_1=2(\alpha+K)$.
  % **only $\alpha+K>1/2$ is necessary**
  Setting $\zeta=\gamma_j \lambda$
  in~(\ref{e:intrepKjp}), we get
 \begin{equation}\label{e:intexprKjp0}
\gamma_j^{p(1-2d)-(2K+1)}\widehat{\kappa}_{j}^{(p)}(\xi)=
\int_{(-\gamma_j\pi,\gamma_j\pi)^p}f_{j}^{(K,p)}(\zeta;\xi)\;\rmd^p\zeta,
\end{equation}
where, for all $j\geq0$, $\lambda\in\mathbb{R}^p$
and $\xi=(\xi_1,\xi_2)\in\mathbb{R}^2$,
\begin{equation*}  %\label{eq:Def-fjkp0}
  f_{j}^{(K,p)}(\gamma_j\lambda;\xi)=\gamma_j^{-2dp-(2K+1)}
  f^{\otimes p}(\lambda)\;
  \widehat{h}_j^{(K)}(\Sigma_p(\lambda)+\xi_1)
  \overline{\widehat{h}_j^{(K)}(\Sigma_p(\lambda)-\xi_2)}\;.
\end{equation*}
Using~(\ref{e:sdf0}),~(\ref{EqLimHj}),~(\ref{e:HjkVSHj}) and $z_{j}\to0$, we
have,  as $j\to\infty$,
  \begin{equation}
    \label{eq:convKappa_0}
f_{j}^{(K,p)}(\zeta;z_{j}/\gamma_j)\to
(f^*(0))^{p}\;
\;\frac{|\widehat{h}_{\infty}(\zeta_1+\cdots+\zeta_{p})|^{2}}
{|\zeta_1+\cdots+\zeta_{p}|^{2K}}\;\prod_{i=1}^p |\zeta_i|^{-2d}\;.
\end{equation}
It turns out, however, that $f_{j}^{(K,p)}(\zeta;z_{j}/\gamma_j)$ cannot be
uniformly bounded by an integrable function over the whole integral domain
$(-\gamma_j\pi,\gamma_j\pi)^p$, but only on a specific subdomain, as we will
show below. By~(\ref{e:sdf0}) and~(\ref{eq:M_replaced_by_K}), we have, for some
constant $C>0$,
\begin{equation}
  \label{eq:bound_fjKp}
\left|f_{j}^{(K,p)}(\zeta;z_{j}/\gamma_j)\right|
\leq C\;\prod_{i=1}^p|\zeta_i|^{-2d}\;
\sup_{|u|\leq |z_j|}
(1+\left|\gamma_j\{(\Sigma_p(\zeta)+u)/\gamma_j\}\right|)^{-2(\alpha+K)}\;.
\end{equation}
The domains are defined using an integer $s$ by taking $\zeta$ such that
$\{(\Sigma_p(\zeta)+u)/\gamma_j\}=(\Sigma_p(\zeta)+u)/\gamma_j-2\pi s$.  In
fact we will use smaller domains that do not depend on $u\in[-|z_j|,|z_j|]$,
namely,
$$
\Gamma_{j}^{(s)}=
\{\zeta\in(-\gamma_j\pi,\gamma_j\pi)^p,\,
-\pi+2\pi s+|z_j|/\gamma_j< \Sigma_p(\zeta)/\gamma_j
< \pi+2\pi s-|z_j|/\gamma_j\}\;.
$$
We note indeed that, for all $\zeta\in\Gamma_{j}^{(s)}$ and
$u\in[-|z_j|,|z_j|]$,
$\{(\Sigma_p(\zeta)+u)/\gamma_j\}=(\Sigma_p(\zeta)+u)/\gamma_j-2\pi s$.  The
following set completes the partition of $(-\gamma_j\pi,\gamma_j\pi)^p$.
$$
\Delta_{j}=\left\{\zeta\in(-\gamma_j\pi,\gamma_j\pi)^p~:~
d\left(\Sigma_p(\zeta)/\gamma_j,\pi+2\pi\mathbb{Z}\right)\leq
|z_j|/\gamma_j\right\}\;,
$$
where $d(x,A)$ denotes the distance between a real $x$ and the set $A$.
We will prove below the following facts.
\begin{enumerate}[(i)]
\item As $j\to\infty$, we have
  \begin{equation}
    \label{eq:convKappa_1}
\int_{\Gamma_{j}^{(0)}}f_{j}^{(K,p)}(\zeta;z_{j}/\gamma_j)\;\rmd\zeta\to
 (f^*(0))^{p}\;\K_p\;.
  \end{equation}
\item If $|s|\geq (p+1)/2$, for $j$ large enough, $\Gamma_{j}^{(s)}$ is an
  empty set.
\item For all $s\neq0$, as $j\to\infty$,
  \begin{equation}
    \label{eq:convKappa_2}
\int_{\Gamma_{j}^{(s)}}f_{j}^{(K,p)}(\zeta;z_{j}/\gamma_j)\;\rmd\zeta\to0\;.
\end{equation}
\item As $j\to\infty$,
  \begin{equation}
    \label{eq:convKappa_3}
\int_{\Delta_{j}}f_{j}^{(K,p)}(\zeta;z_{j}/\gamma_j)\;\rmd\zeta\to0\;.
\end{equation}
\end{enumerate}
To conclude the proof, we show (i), (ii), (iii) and (iv) successively.

First consider (i). It follows from~(\ref{eq:bound_fjKp}), the
definition of $\Gamma_{j}^{(0)}$ and $|z_j|\to0$ that, for $j$ large enough,
$$
\1_{\Gamma_{j}^{(0)}}(\zeta)\left|f_{j}^{(K,p)}(\zeta;z_{j}/\gamma_j)\right|
\leq C\;\prod_{i=1}^p|\zeta_i|^{-2d}\;(1/2+\left|\Sigma_p(\zeta)\right|)^{-2(\alpha+K)}\;.
$$
Observe that, by Lemma~\ref{lem:lem8.3CRTT}, and since $\alpha>1/2$, $K\geq0$
and $p(1-2d)<1$, the right-hand side of the last display is integrable.
Then~(\ref{eq:convKappa_1}) follows from~(\ref{eq:convKappa_0})
and the dominated convergence theorem.

Assertion (ii) follows from the definition of $\Gamma_{j}^{(s)}$.

We now prove (iii) and thus take $s\neq0$. Using~(\ref{eq:bound_fjKp}) and
$|z_j|\to0$, we get, for all $\zeta\in\Gamma_{j}^{(s)}$ and $j$ large enough,
$$
\left|f_{j}^{(K,p)}(\zeta;z_{j}/\gamma_j)\right|
\leq C\;\prod_{i=1}^p|\zeta_i|^{-2d}\;
(1/2+\left|\Sigma_p(\zeta)-2\pi s\gamma_j\right|)^{-2(\alpha+K)}\;.
$$
The limit~(\ref{eq:convKappa_2}) then follows from Lemma~\ref{lem:lem8.4CRTT}
applied with $q=p$, $M_1=2(K+\alpha)$ and
$a=2\pi\gamma_j s$.

Finally we prove Assertion~(iv). In this case, we observe
that~(\ref{eq:bound_fjKp}) and implies
$$
|f_{j}^{(K,p)}(\zeta;z_{j}/\gamma_j)|\leq C\prod_{i=1}^p|\zeta_i|^{-2d}\;.
$$
This bound and Lemma~\ref{lem:lem8.3CRTT} yields
$$
\int_{\Delta_{j}}f_{j}^{(K,p)}(\zeta;z_{j}/\gamma_j)\;\rmd\zeta
\leq
C\;\int_{-p\gamma_j\pi}^{p\gamma_j\pi}
\1_{d(t/\gamma_j,\pi+2\pi\mathbb{Z})\leq |z_j|/\gamma_j}\;\rmd t=O(|z_j|)\;.
$$
Hence, since $|z_j|\to0$, we obtain~(\ref{eq:convKappa_3}) and the proof is
achieved.
\end{proof}

\subsection{Other technical lemmas}
\begin{lemma}\label{lem:Dn}
Define the Dirichlet kernel $D_n$ as in~(\ref{eq:dirichlet}). Then
\begin{equation}\label{eq:DnBound}
\sup_{\theta\in \mathbb{R}} \sup_{n\geq1}
\; (1+|n\{\theta/n\}|)\left|D_n(\theta/n)\right| < \infty \;.
\end{equation}
\end{lemma}
\begin{proof}
We observe that $|\rme^{\rmi\lambda}-1|\geq 2|\{\lambda\}|/\pi$.
Hence,  for all $\theta\in \mathbb{R}$,
\begin{equation*}%\label{eq:DnBound}
\left|D_n(\theta/n)\right|
\leq
\frac{\pi}{2}
\frac{|\rme^{\rmi\theta}-1|}{|n\{\theta/n\}|}
=\frac{\pi}{2}
\frac{|\rme^{\rmi n \{\theta/n\} }-1|}{|n\{\theta/n\}|} \;.
\end{equation*}
(We use the usual continuous extension convention $(\rme^{\rmi 0}-1)/0=1$).
Now, using that $|e^{\rmi u}-1|\leq 2|u|/(1+|u|)$ on $u\in\mathbb{R}$, we
get~(\ref{eq:DnBound}).
\end{proof}

\begin{lemma}\label{lem:lem8.3CRTT}
Let $p$ be a positive integer and $f:\mathbb{R}\to\mathbb{R}_+$.
Then, for any $\beta\in\mathbb{R}^{q}$,
  \begin{equation}
    \label{eq:chagtVar}
\int_{\mathbb{R}^q}f(y_1+\cdots+y_q)\;
\prod_{i=1}^{q}|y_i|^{\beta_i}\; \rmd y_1\cdots\rmd
y_q=\Gamma\times \int_{\mathbb{R}}
f(s)|s|^{q-1+\beta_1+\cdots+\beta_q} \rmd s\;,
  \end{equation}
where, for all $i\in\{1,\cdots,q\}$,
$B_i=\beta_i+\cdots+\beta_{q}$ and
$$
\Gamma=\prod_{i=2}^{q}
\left(\int_{\mathbb{R}}|t|^{q-i+B_i}|1-t|^{\beta_{i-1}}\rmd
t\right)\;.
$$
(We note that $\Gamma$ may be infinite in which
case~(\ref{eq:chagtVar}) holds with the convention
$\infty\times0=0$).
\end{lemma}
\begin{proof}
This follows from Lemma~8.3 in~\cite{clausel-roueff-taqqu-tudor-2011a}.
\end{proof}

\begin{lemma}\label{lem:lem8.4CRTT}
  Let $d\in(0,1/2)$ and $q$ be a positive integer such that $q<1/(1-2d)$.
  Let $M_1>1$. Set for any $a\in\mathbb{R}$,
\[
J_q(a;M_1;d)=\int_{\mathbb{R}^{q}}
(1+|\Sigma_q(\zeta)-a|)^{-M_1}\prod\limits_{i=1}^{q}|\zeta_i|^{-2d}
\;\rmd\zeta.
\]
Then one has
\begin{equation}\label{e:Jq1}
\sup_{a\in\mathbb{R}}(1+|a|)^{1-q(1-2d)}J_q(a;M_1;d)<\infty\;.
\end{equation}
In particular,
\[
J_q(0;M_1;d)<\infty,
\]
and
\[
J_q(a;M_1;d)=O(|a|^{-(1-q(1-2d)})\quad\mbox{ as }a\to\infty\;.
\]
\end{lemma}
\begin{proof}
This follows from Lemma~8.4 of in~\cite{clausel-roueff-taqqu-tudor-2011a}.
\end{proof}

\begin{lemma}\label{lem:Ja}
 Define, for all $a>0$ and $\beta_1\in(0,1)$,
\begin{equation}\label{e:J1a}
J_{1,a}(s_1;\beta_1)=|s_1|^{-\beta_1},\quad s_1\in\mathbb{R}\;,
\end{equation}
and, for any integer $m\geq2$ and  $\beta=(\beta_1,\cdots,\beta_m)\in
  (0,1)^m$,
\begin{equation}\label{e:Jma}
J_{m,a}(s_1;\beta)=\int_{s_2=-(m-1)a}^{(m-1)a}\dots\int_{s_m=-a}^{a}
\prod_{i=2}^m|s_{i-1}-s_i|^{-\beta_{i-1}}\; |s_m|^{-\beta_m}
\;\rmd s_m\dots \rmd s_2,\quad s_1\in\mathbb{R}\;.
\end{equation}
Then
\begin{enumerate}[(i)]
\item if $\beta_1+\dots+\beta_m> m-1$, one has
$$
C_m(\beta)=
\sup_{a>0}\;\sup_{s_1\in\mathbb{R}}\left(  |s_1|^{-(m-1-(\beta_1+\dots+\beta_m))}
J_{m,a}(s_1;\beta) \right) < \infty \;,
$$
\item if $\beta_1+\dots+\beta_m=m-1$, one has
$$
C_m(\beta)=
\sup_{a>0}\;\sup_{|s_1|\leq ma}\left( \frac1{1+\log(ma/|s_1|)}
J_{m,a}(s_1;\beta) \right)< \infty \;,
$$
\item if there exists $q\in\{2,\dots,m\}$ such that
  $\beta_q+\dots+\beta_m=m-q$, one has
$$
C_m(\beta)=
\sup_{a>0}\;\sup_{|s_1|\leq ma} \left(
  \frac{a^{-(q-1-(\beta_{1}+\dots+\beta_{q-1}))}}{1+\log(ma/|s_1|)}J_{m,a}(s_1;\beta)  \right)< \infty \;,
$$
\item if  $\beta_1+\dots+\beta_m<m-1$ and for all $q\in\{1,\dots,m-1\}$, we
  have  $\beta_q+\dots+\beta_m\neq m-q$, one has
$$
C_m(\beta)=
\sup_{a>0}\;\sup_{|s_1|\leq ma} \left(a^{-(m-1-(\beta_1+\dots+\beta_m))}J_{m,a}(s_1;\beta)  \right)< \infty \;.
$$
\end{enumerate}
\end{lemma}
\begin{remark}
  We observe that Cases (ii),(iii) and (iv) can be put together as the
  following formula, valid for all $\beta\in(0,1)^m$ such that
  $\beta_1+\dots+\beta_m\leq m-1$,
  \begin{equation}
    \label{eq:Jacases234}
    C_m(\beta)=
    \sup_{a>0}\;\sup_{|s_1|\leq ma}
    \left(
      \frac{a^{-(q-1-(\beta_{1}+\dots+\beta_{q-1}))}}{\{1+\log(ma/|s_1|)\}^\varepsilon}
      J_{m,a}(s_1;\beta)  \right)< \infty \;,
  \end{equation}
where $\varepsilon=1$ if there exists $q\in\{1,\dots,m\}$ such that
  $\beta_q+\dots+\beta_m=m-q$, and  $\varepsilon=0$ otherwise.
  We may also include case (i) as follows,
  \begin{equation}
    \label{eq:Jacases1234}
    C_m(\beta)=
    \sup_{a>0}\;\sup_{|s_1|\leq ma}
    \left(
      \frac{a^{-(m-1-(\beta_{1}+\dots+\beta_{m}))_+}
|s_1|^{(m-1-(\beta_{1}+\dots+\beta_{m}))_-}}
{\{1+\log(ma/|s_1|)\}^\varepsilon}
      J_{m,a}(s_1;\beta)  \right)< \infty \;,
  \end{equation}
where $\varepsilon$ is as above, and $a_+=\max(a,0)$ and $a_-=\max(-a,0)$ denote
the positive and negative parts of $a$, respectively.
\end{remark}
\begin{proof}
Observe first that for all $m\geq1$,
  \begin{equation}
    \label{eq:inductionJ}
J_{m,a}(s_1;\beta)=\int_{s_2=-(m-1)a}^{(m-1)a}|s_2-s_1|^{-\beta_1}\;
J_{m-1,a}(s_2;\beta')\;\rmd s_2\;,
  \end{equation}
  where $\beta'=(\beta_2,\dots,\beta_{m})$. The bounds $C_m(\beta)$ in the different cases
  will follow by induction on $m$.

  Let us first prove the result for $m=1$ and $m=2$. If $m=1$,
  $\beta=\beta_1\in(0,1)$ only satisfies the condition of Case (i) and, since
  $J_{1,a}$ is given by~(\ref{e:J1a}), the result holds for
  $m=1$. Assume now that $m=2$ and $s_1\neq 0$ and set $s_2=v |s_1|$. Then
\begin{equation}\label{eq:chgtvarJa}
J_{2,a}(s_1;\beta)= |s_1|^{1-(\beta_1+\beta_2)}\int_{-a/|s_1|}^{a/|s_1|}\frac{\rmd v}{|1-v|^{\beta_1}|v|^{\beta_2}}\;.
\end{equation}
In the case $\beta_1+\beta_2>1$, we are in Case~(i). Since
$\int_{\mathbb{R}}\frac{\rmd v}{|1-v|^{\beta_1}|v|^{\beta_2}}$ is finite, the
required upper bound holds. If $\beta_1+\beta_2\leq 1$, we are either in
Case~(ii) or (iv) and the result follows from the following bounds valid for
some constant $c$ depending only on $\beta$, if $\beta_1+\beta_2< 1$ and
$x\geq1/2$,
$$
\int_{-x}^{x}\frac{\rmd v}{|1-v|^{\beta_1}|v|^{\beta_2}}\leq
c x^{1-(\beta_1+\beta_2)}\;,
$$
and, if $\beta_1+\beta_2= 1$ and $x\geq1/2$,
$$
\int_{-x}^{x}\frac{\rmd v}{|1-v|^{\beta_1}|v|^{\beta_2}}\leq C
(1+\log(2x))\;.
$$
This prove the result for $m=2$ because $x=a/|s_1|\geq 1/2$.

Let us now assume that the result holds for some positive integer
$m-1$ and prove it for $m$. We consider two different cases.
\begin{enumerate}
\item If $\beta$ satisfies the conditions of Case (i), Case (ii), or Case (iv)
  then $\beta'$ satisfies the conditions of Case (i) or (iv). Then
  by~(\ref{eq:inductionJ}) and the induction
  assumption,
\[
J_{m,a}(s_1;\beta)\leq C_{m-1}(\beta')
a^{[m-2-\Sigma_{m-1}(\beta')]_+}
\int_{-(m-1)a}^{(m-1)a}
|s_2-s_1|^{-\beta_1}|s_2|^{-[\Sigma_{m-1}(\beta')-(m-2))]_+}\rmd s_2\;,
\]
where $\Sigma_{m-1}(\beta')=\beta_2+\dots+\beta_m$ and $[x]_+=\max(x,0)$.  If
$\Sigma_{m-1}(\beta')<m-2$ (so that $\beta$ satisfies (iv)), the conclusion
follows from the following bound
valid for some constant $c$ depending only on $\beta$ and all $x\geq|s_1|/2$,
$$
\int_{-x}^{x}
|s_2-s_1|^{-\beta_1}\rmd s_2
=|s_1|^{1-\beta_1}\int_{-x/|s_1|}^{x/|s_1|}|u-1|^{-\beta_1}\rmd u
\leq c x^{1-\beta_1}\;.
$$
Now if $\Sigma_{m-1}(\beta')>m-2$, we observe that
\[
\int_{-(m-1)a}^{(m-1)a}|s_2-s_1|^{-\beta_1}|s_2|^{-[\beta_2+\cdots+\beta_m-(m-2)]}\rmd s_2=J_{2,(m-1)a}(s_1;\beta_1,\beta_2+\cdots+\beta_m-(m-2))\;.
\]
The upper bound of $J_{m,a}(s_1;\beta)$ then follows from the case $m=2$.
\item If $\beta$ satisfies the condition of Case (iii), then $\beta'$ either satisfies the conditions of Case (ii) or (iii). The proof is exactly similar to this just above up to a logarithmic correction.
\end{enumerate}

\end{proof}

\begin{lemma}\label{lem:Jnbound}
Let $S>1$ and $(\beta_1,\beta_2)\in[0,1)^2$ such that $\beta_1+\beta_2<1$, and set
$g_i(t)=|t|^{-\beta_i} (1+|t|)^{\beta_i-S}$. Then
\begin{equation}
  \label{eq:JnBound}
\sup_{\nu\geq0}
\left(\nu
\int_{\mathbb{R}^2}(1+\nu|\{w_1+w_2\}|)^{-2}g_1(w_1)g_2(w_2)\;\rmd w\right)
< \infty \;.
\end{equation}
\end{lemma}
\begin{proof}
  Denote by $J(\nu)$ the quantity in parentheses in~(\ref{eq:JnBound}).  We
  denote here by $C$ a positive constant that may change from line to line, but
  whose value does not depend on $\nu$. Setting $u=w_1+w_2$ in the integral
  with respect to $w_1$ and then integrating with respect to $w_2$, Lemma 8.1
  in~\cite{clausel-roueff-taqqu-tudor-2011a} yields
$$
J(\nu)\leq C\nu\int_{u\in\mathbb{R}}(1+\nu|\{u\}|)^{-2}\,(1+|u|)^{-S}\rmd u\;.
$$
Since the integral is bounded independently of $\nu$, $J$ is bounded on compact
subsets of $[0,\infty)$, hence we may consider $\nu\geq2$ in the remainder of
the proof.  We shall use the bound $1+x\geq \max(1,x)$ for $x\geq0$. Splitting
the integral of the last display on the two domains defined by the position of
$|\{u\}|$ with respect to $\nu^{-1}$, we get $J(\nu)\leq C(J_1(\nu)+J_2(\nu))$,
with
$$
J_1(\nu)=\nu\int_{|\{u\}|\leq\nu^{-1}} (1+|u|)^{-S} \;\rmd u\;,
$$
and
$$
J_2(\nu)= \nu^{-1}\int_{|\{u\}|\geq\nu^{-1}}
|\{u\}|^{-2}\;(1+|u|)^{-S}\;\rmd u\;.
$$
We have
$$
J_1(\nu)=\nu\sum_{k\in\mathbb{Z}}
\int_{2k\pi-\nu^{-1}}^{2k\pi+\nu^{-1}}(1+|u|)^{-S} \rmd u\;.
$$
For $\nu\geq2$ the integral in the parentheses of the last display is
less than
$2\nu^{-1}(1/2+|2k\pi|)^{-S}$. Since $S>1$,
we get that $J_1(u)$ is bounded over the domain $\nu\geq2$.

It remains to prove that $J_2(\nu)$ is bounded for $\nu$ large enough.
We have, setting $v=u-2k\pi$ for each $k$,
$$
J_2(\nu)= \nu^{-1}\sum_{k\in\mathbb{Z}}
\int_{\nu^{-1}\leq |v|\leq \pi}|v|^{-2}(1+|2k\pi+v|)^{-S}\;\rmd v\;.
$$
Now since
$$
\sup_{v\in\mathbb{R}}\sum_{k\in\mathbb{Z}}(1/2+|2k\pi+v|)^{-S}<\infty\;,
$$
we get by inverting the integral with the summation,
$$
J_2(\nu)\leq C \nu^{-1}\int_{\nu^{-1}\leq |v|\leq \pi}|v|^{-2}\;\rmd v\;.
$$
Hence $J_2$ is bounded over the domain $\nu\geq2$, completing the
proof.
\end{proof}

\section{Integral representations}\label{s:appendix}
It is convenient to use an integral representation in the spectral
domain to represent the random processes (see for
example~\cite{major:1984,nualart:2006}). The stationary Gaussian
process $\{X_k,k\in\mathbb{Z}\}$ with spectral
density~(\ref{e:sdf}) can be written as
\begin{equation}\label{e:intrepX}
X_\ell=\int_{-\pi}^{\pi}\rme^{\rmi\lambda
  \ell}f^{1/2}(\lambda)\rmd\widehat{W}(\lambda)=\int_{-\pi}^{\pi}\frac{\rme^{\rmi\lambda
    \ell}f^{*1/2}(\lambda)}{|1-\rme^{-{\rmi}\lambda}|^d}\rmd\widehat{W}(\lambda),\quad\ell\in\mathbb{Z}\;.
\end{equation}
This is a special case of
\begin{equation}\label{e:int}
\widehat{I}(g)=\int_{\mathbb{R}}g(x)\rmd\widehat{W}(x),
\end{equation}
where $\widehat{W}(\cdot)$ is a complex--valued Gaussian random measure
satisfying, for any Borel sets $A$ and $B$ in $\mathbb{R}$,
$\mathbb{E}(\widehat{W}(A))=0$,
$\mathbb{E}(\widehat{W}(A)\overline{\widehat{W}(B)})=|A\cap B|$ and
$\widehat{W}(A)=\overline{\widehat{W}(-A)}$.
The integral~(\ref{e:int}) is defined for any function $g\in L^2(\mathbb{R})$
and one has the isometry
\[
\mathbb{E}(|\widehat{I}(g)|^2)=\int_{\mathbb{R}}|g(x)|^2\rmd x\;.
\]
The integral $\widehat{I}(g)$, moreover, is real--valued if
$g(x)=\overline{g(-x)}$.

We shall also consider multiple It\^{o}--Wiener integrals
\[
\widehat{I}_q(g)=\int^{''}_{\mathbb{R}^q}g(\lambda_1,\cdots,\lambda_q)\rmd \widehat{W}(\lambda_1)\cdots\rmd\widehat{W}(\lambda_q)
\]
where the double prime indicates that one does not integrate on hyperdiagonals
$\lambda_i=\pm \lambda_j,i\neq j$. The integrals $\widehat{I}_q(g)$ are handy
because we will be able to expand our non--linear functions $G(X_k)$ introduced
in Section~\ref{s:intro} in multiple integrals of this type.

These multiples integrals are defined for
$g\in\overline{L^2}(\mathbb{R}^{q},\mathbb{C})$, the space of complex valued
functions defined on $\mathbb{R}^{q}$ satisfying
\begin{gather}\label{e:antisym}
g(-x_{1},\cdots,-x_{q})= \overline{g(x_{1},\cdots, x_{q})}\mbox{ for }(x_{1},\cdots, x_{q}) \in \mathbb{R}^q\;,\\
\label{e:fL2}
\Vert g\Vert ^{2} _{L^{2}}:= \int_{\mathbb{R}^{q}} \left| g(x_{1},\cdots, x_{q}) \right| ^{2} \rmd x_{1}\cdots \rmd x_{q} <\infty\;.
\end{gather}
The integral $\widehat{I}_q(g)$ is real valued and verifies
$\widehat{I}_q(g)=\widehat{I}_q(\tilde{g})$,
where
$$
\tilde{g} (x_{1},\cdots, x_{q}) =\frac1{q!}\sum_{\sigma}
g(x_{\sigma (1)},\cdots, x_{\sigma (q)})\;.
$$
Here the sum is over all permutations of $\{1,\dots,q\}$.
\begin{equation}\label{e:cov-multiple}
\mathbb{E}(\widehat{I}_{q}(g_1) \widehat{I}_{q'}(g_2))=
\begin{cases}
q! \langle \tilde{g_1}, \tilde{g_2}\rangle _{L^{2}}&\text{if } q=q'\\
0& \text{if } q\neq q'.
\end{cases}
\end{equation}
Hermite polynomials are related to multiple integrals as follows : if
$X=\int_{\mathbb{R}}g(x)\rmd\widehat{W}(x)$ with
$\mathbb{E}(X^2)=\int_{\mathbb{R}}|g(x)|^2\rmd x=1$ and $g(x)=\overline{g(-x)}$
so that $X$ has unit variance and is real--valued, then
\begin{equation}\label{e:herm-integ}
H_q(X)=\widehat{I}_q(g^{\otimes q})=\int_{\mathbb{R}^q}^{''}g(x_1)\cdots g(x_q)\rmd \widehat{W}(x_1)\cdots\rmd\widehat{W}(x_q)\;.
\end{equation}
Since $X$ has unit variance, one has for any $\ell\in\mathbb{Z}$,
\begin{align*}
H_q(X_\ell)&=
H_q\left(\int_{-\pi}^{\pi}\rme^{\rmi\xi \ell}f^{1/2}(\xi)
\rmd\widehat{W}(\xi)\right)\\
&=\int_{(-\pi,\pi]^q}^{''}
\rme^{\rmi\ell(\xi_1+\cdots+\xi_q)}\times
\left(f^{1/2}(\xi_1)\times\cdots\times f^{1/2}(\xi_q)\right)
\;\rmd\widehat{W}(\xi_1)\cdots \rmd\widehat{W}(\xi_q)\;.
\end{align*}
Then by~(\ref{e:W2}), we have
\begin{equation}\label{eq:Wjk_representation}
W_{j,k}=\sum_{\ell\in\mathbb{Z}}h_j^{(K)}(\gamma_j
k-\ell)H_{q_0}(X_\ell)=\widehat{I}_{q_0}(f_{j,k}^{(q_0)})
\end{equation}
with
\begin{equation}\label{e:fjk_representation1}
f_{j,k}^{(q)}(\xi_1,\cdots,\xi_q)= \rme^{\rmi
k\gamma_j(\xi_1+\cdots+\xi_q)}\times
\widehat{h}_{j}^{(K)}(\xi_1+\cdots+\xi_q) f^{1/2}(\xi_1)\cdots
f^{1/2}(\xi_q) \1_{(-\pi,\pi)}^{\otimes q}(\xi) \;,
\end{equation}
because by~(\ref{e:dF}),
\begin{eqnarray*}
\sum_{\ell\in\mathbb{Z}}\rme^{\rmi\ell(\xi_1+\cdots+\xi_q)}h_j^{(K)}(\gamma_j k-\ell)&=&\rme^{\rmi\gamma_j
k(\xi_1+\cdots+\xi_q)}\sum_{u\in\mathbb{Z}}\rme^{-\rmi
u(\xi_1+\cdots+\xi_q)}h_j^{(K)}(u)\\
&=&\rme^{\rmi\gamma_j
k(\xi_1+\cdots+\xi_q)}\widehat{h}_j^{(K)}(\xi_1+\cdots+\xi_q)\;.
\end{eqnarray*}

Observe now that since we have defined the Fourier transform of a function $f\in L^2(\mathbb{R}^q)$ as
\[
\widehat{f}(\xi)=\int_{\mathbb{R}^q}f(x)\rme^{-\rmi x\xi}\rmd x\in \overline{L^2}(\mathbb{R}^q)\;,
\]
we have by Parseval
\[
\|\widehat{f}\|^2_{\overline{L^2}(\mathbb{R}^q)}=(2\pi)^q\;\|f\|^2_{L^2(\mathbb{R}^q)}\;.
\]

Since moreover, $\mathbb{E}(I_q(\widehat{f})^2)=\|\widehat{f}\|^2_{\overline{L^2}(\mathbb{R}^q)}$ and $\mathbb{E}(\widehat{I}_q(f)^2)=\|f\|^2_{L^2(\mathbb{R}^q)}$, we have
\begin{equation}\label{e:four}
I_q(\widehat{f})\overset{(\mathcal{L})}{=}(2\pi)^{q/2}\widehat{I}_q(f)\;.
\end{equation}

The following proposition is an
extension to our complex--valued setting of a corresponding result
in~\cite{nualart:2006} for multiple integrals in a real--valued setting. Since
it plays an essential role, we provide a proof for the convenience of the
reader.
\begin{proposition}\label{pro:intproduct}
  Let $(q,q')\in \mathbb{N}^2$. Assume that $f,g$ are two symmetric functions
  belonging respectively to $\overline{L^2}(\mathbb{R}^q)$ and
  $\overline{L^2}(\mathbb{R}^{q'})$ then the following product formula holds :
\begin{equation}\label{EqProdIntSto}
\widehat{I}_q(f)\widehat{I}_{q'}(g)=\sum\limits_{p=0}^{q\wedge q'}
 p
!\begin{pmatrix}q\\p\end{pmatrix}\begin{pmatrix}q'\\p\end{pmatrix}\widehat{I}_{q+q'-2p}(f\overline{\otimes}_p
g),
\end{equation}
where $f\overline{\otimes}_0 g=f\otimes g$ is the usual tensor product and, for
any $p\in\{1,\cdots,q\wedge q'\}$,
\begin{equation}\label{e:times-p}
(f\overline{\otimes}_p
g)(t_1,\cdots,t_{q+q'-2p})=\int_{\mathbb{R}^p}f(t_1,\cdots,t_{q-p},s)g(t_{q-p+1},\cdots,t_{q+q'-2p},-s)\rmd^p
s\;.
\end{equation}
\end{proposition}
\begin{proof}
We first assume that $f$ and $g$ are of the form
\[
f=f_1\otimes f_2,\,g=g_1\otimes g_2\;,
\]
where $f_1,f_2,g_1,g_2$ belong respectively to
$\overline{L^2}(\mathbb{R}^{q-p},\mathbb{C}),\overline{L^2}(\mathbb{R}^{q'-p},\mathbb{C}),\overline{L^2}(\mathbb{R}^{p},\mathbb{C}),\overline{L^2}(\mathbb{R}^{p},\mathbb{C})$. In
that special case, using that for any $q\geq 1$ and any $f\in
\overline{L^2}(\mathbb{R}^{q})$,
$\widehat{I}_q(f)=(2\pi)^{-q/2}I_q(\widehat{f})$ by~(\ref{e:four}), one has
\begin{equation}\label{e:f1}
\widehat{I}_q(f)\widehat{I}_{q'}(g)=\widehat{I}_q(f_1\otimes
f_2)\widehat{I}_{q'}(g_1\otimes g_2)
=(2\pi)^{-(q+q')/2}I_q(\widehat{f_1}\otimes \widehat{f_2})I_{q'}(\widehat{g_1}\otimes \widehat{g_2})\;.
\end{equation}
The assumptions on functions $f_1,f_2,g_1,g_2$ imply that their Fourier
transform $\widehat{f_1}, \widehat{f_2}, \widehat{g_1}, \widehat{g_2}$ are
real--valued functions belonging respectively to
$L^2(\mathbb{R}^{q-p},\mathbb{R}),\,L^2(\mathbb{R}^{q'-p},\mathbb{R}),\,L^2(\mathbb{R}^{p},\mathbb{R})$
and $L^2(\mathbb{R}^{\ell},\mathbb{R})$. Then one can apply the usual product
formula for multiple Wiener-It\^o integrals (see for example~\cite{nualart:2006}) and
deduce that :
\begin{equation}\label{e:f2}
I_q(\widehat{f_1}\otimes \widehat{f_2})I_{q'}(\widehat{g_1}\otimes \widehat{g_2})=\sum\limits_{p=0}^{q\wedge q'}p !\begin{pmatrix}q\\p\end{pmatrix}\begin{pmatrix}q'\\p\end{pmatrix}I_{q+q'-2p}((\widehat{f_1}\otimes \widehat{f_2})\otimes_p (\widehat{g_1}\otimes \widehat{g_2}))\;.
\end{equation}
Note now that for any $p$
\begin{align*}
(\widehat{f_1}\otimes \widehat{f_2})\otimes_p (\widehat{g_1}\otimes \widehat{g_2})&=
\int_{\mathbb{R}^p}\widehat{f_1}(t_1,\cdots,t_{q-p})\widehat{f_2}(s)\widehat{g_1}(t_{q-p+1},\cdots,t_{q+q'-2p})\widehat{g_2}(s)\rmd s\\
&=\widehat{f_1}(t_1,\cdots,t_{q-p})\widehat{g_1}(t_{q-p+1},\cdots,t_{q+q'-2p})
\int_{\mathbb{R}^p}\widehat{f_2}(s)\widehat{g_2}(s)\rmd s\\
&=\widehat{f_1}(t_1,\cdots,t_{q-p})\widehat{g_1}(t_{q-p+1},\cdots,t_{q+q'-2p})
(2\pi)^p \int_{\mathbb{R}^p}f_2(t)\overline{g_2(t)}\rmd t\\
&=\widehat{f_1}(t_1,\cdots,t_{q-p})\widehat{g_1}(t_{q-p+1},\cdots,t_{q+q'-2p})
(2\pi)^p \int_{\mathbb{R}^p}f_2(t)g_2(-t)\rmd t  \;,
\end{align*}
since $\overline{g_2(t)}=g_2(-t)$ and using the Parseval's formula.
Hence
\begin{align*}
I_{q+q'-2p}((\widehat{f_1}\otimes \widehat{f_2})\otimes_p (\widehat{g_1}\otimes \widehat{g_2}))
&=(2\pi)^p\left(\int_{\mathbb{R}^p}f_2(t)g_2(-t)dt\right)\times I_{q+q'-2p}(\widehat{f_1}\otimes\widehat{g_1})\\
&=(2\pi)^p\left(\int_{\mathbb{R}^p}f_2(t)g_2(-t)dt\right)\times I_{q+q'-2p}(\widehat{f_1\otimes g_1})\\
&=(2\pi)^p\left(\int_{\mathbb{R}^p}f_2(t)g_2(-t)dt\right)\times (2\pi)^{(q+q'-2p)/2}\widehat{I}_{q+q'-2p}(f_1\otimes g_1)\\
&=(2\pi)^{(q+q')/2}\left(\int_{\mathbb{R}^p}f_2(t)g_2(-t)dt\right)\widehat{I}_{q+q'-2p}(f_1\otimes g_1)\\
&=(2\pi)^{(q+q')/2}\widehat{I}_{q+q'-2p}(f\overline{\otimes}_p g)\;.
\end{align*}
Using the last equality and equations~(\ref{e:f1}), (\ref{e:f2}),we get the claimed results for this special case.  The conclusion for
general
$f$ and $g$ follows using the density of
$L^2(\mathbb{R}^{q-p},\mathbb{R})\otimes L^2(\mathbb{R}^{p},\mathbb{R})$ in
$L^2(\mathbb{R}^q,\mathbb{R})$.
\end{proof}

\bibliographystyle{plainnat}
\bibliography{lrd}
\end{document}